\documentclass{amsart}
\usepackage{amssymb,latexsym,epic} 
\newlength{\fw}
\newlength{\hs} 
\newlength{\agk}
\newlength{\agku}
\newlength{\jmr}
\newlength{\bernd}
\newtheorem{lemma}{Lemma}
\newtheorem{cor}{Corollary}
\newtheorem{prop}{Proposition}
\newtheorem{propdfn}[prop]{Proposition and Definition}

\newtheorem{main}{Main Theorem}

\newtheorem{mres}{Main Result}
 
\newtheorem{anti}{Antipodality Theorem}
 
\newtheorem{alg}{Algorithm}
\newtheorem{thm}{Theorem}

\newtheorem{afpii}{Affine Point Theorem II}

\newtheorem{rem}{Remark}

\newtheorem{ex}{Example}
\newtheorem{dfn}[prop]{Definition}

\newcommand{\gen}{\mathrm{gen}}
\newcommand{\sm}{\mathrm{sum}}
\newcommand{\eps}{\varepsilon}
\newcommand{\cA}{\mathcal{A}}
\newcommand{\cD}{\mathcal{D}}
\newcommand{\wcD}{\widetilde{\cD}}
\newcommand{\Ds}{\cD_{\mathrm{shift}}}

\newcommand{\cE}{\mathcal{E}}
\newcommand{\cH}{\mathcal{H}}

\newcommand{\cN}{\mathcal{N}}
\newcommand{\cO}{\mathcal{O}}

\newcommand{\cI}{\mathcal{I}}
\newcommand{\cP}{\mathcal{P}}
\newcommand{\cV}{\mathcal{V}}
\newcommand{\cX}{\mathcal{X}}
\newcommand{\cY}{\mathcal{Y}}
\newcommand{\V}{\mathrm{Vol}}
\newcommand{\supp}{\mathrm{Supp}}
\newcommand{\conv}{\mathrm{Conv}}
\newcommand{\init}{\mathrm{in}}
\newcommand{\thth}{{\underline{\mathrm{th}}}} 
\newcommand{\st}{{\underline{\mathrm{st}}}} 
\newcommand{\nd}{{\underline{\mathrm{nd}}}} 
\newcommand{\Pro}{\mathbb{P}}
\newcommand{\Q}{\mathbb{Q}}
\newcommand{\R}{\mathbb{R}}
\newcommand{\C}{\mathbb{C}}
\newcommand{\G}{\mathbb{G}}
\newcommand{\bL}{\mathbb{L}}
\newcommand{\N}{\mathbb{N}}
\newcommand{\No}{\N\cup\{0\}}
\newcommand{\Non}{(\No)^n}
\newcommand{\Z}{\mathbb{Z}}

\newcommand{\ch}{\mathrm{char}}

\newcommand{\ord}{\mathrm{ord}}
\newcommand{\area}{\mathrm{Area}}
\newcommand{\spec}{\mathrm{Spec}}
\newcommand{\divisor}{\mathrm{Div}}

\newcommand{\chow}{\mathrm{Chow}}
\newcommand{\relint}{\mathrm{RelInt}}
\newcommand{\fan}{\mathrm{Fan}}
\newcommand{\hyp}{\mathrm{Hyper}}
\newcommand{\lin}{\mathrm{Lin}}
\newcommand{\res}{\mathrm{Res}}
\newcommand{\newt}{\mathrm{Newt}}

\newcommand{\Pc}{{P^\llcorner}} 
\newcommand{\wP}{{\widetilde{P}}} 
\newcommand{\Sn}{\mathcal{S}^{n-1}}
\newcommand{\Zn}{\Z^n}
\newcommand{\Qn}{\Q^n}
\newcommand{\Rn}{\R^n}

\newcommand{\Rni}{\Rn\!\setminus\!\sigma_I} 
 
\newcommand{\gln}{\G\bL_n}
\newcommand{\glnz}{\gln(\Z)}
\newcommand{\Cn}{\C^n}

\newcommand{\Ks}{K^*}
\newcommand{\Kn}{K^n}
\newcommand{\Kni}{\Kn\!\setminus\!\hyp(I)}

\newcommand{\Knjp}{\Kn\!\setminus\!\hyp(J')} 
\newcommand{\Cni}{\Cn\!\setminus\!\hyp(I)} \newcommand{\bO}{{\bf O}}
\newcommand{\Qns}{\Qn\!\setminus\!\{\bO\}} 
\newcommand{\Rns}{\Rn\!\setminus\!\{\bO\}} \newcommand{\Csn}{{(\C^*)}^n}
\newcommand{\Ksn}{{(K^*)}^n}
 
\newcommand{\cM}{\mathcal{M}}
\newcommand{\cS}{\mathcal{S}}
\newcommand{\cSM}{\cS\cM}

\newcommand{\cU}{\mathcal{U}}
\newcommand{\cB}{\mathcal{B}}
\newcommand{\cC}{\mathcal{C}}
\newcommand{\cF}{\mathcal{F}}

\newcommand{\cT}{\mathcal{T}}
\newcommand{\wcT}{\widetilde{\cT}}

\begin{document}

\setlength{\oddsidemargin}{.25in}
\setlength{\evensidemargin}{.25in}

\title[Toric Intersection Theory for Affine Root Counting]{Toric
Intersection Theory for Affine Root Counting$^1$}
\footnotetext[1]{Submitted~to~the~Journal~of~Pure~and~Applied~Algebra.~Also~ 
\mbox{available}~\mbox{on-line}~at~{\tt http://www-math.mit.edu/\~{}rojas}.}

\author{J. Maurice Rojas}
\thanks{This research was completed at MSRI and partially funded by a
Chateaubriand Fellowship, an NSF Mathematical Sciences Postdoctoral 
Research Fellowship, and NSF grant DMS-9022140.}

\address{Massachusetts Institute of Technology\\
Mathematics Department\\
77 Mass.\ Ave.\\
Cambridge, MA \ 02139, U.S.A. }

\email{rojas@math.mit.edu} 

%
%
%

\date{October 27, 1996} 

\dedicatory{This paper is dedicated to Basil Gordon on the occassion of his
$2^{6\thth}$ birthday.}

\begin{abstract} 
Given any polynomial system with fixed monomial term structure, we give 
explicit formulae for the generic number of roots with specified 
coordinate vanishing restrictions. For the case of affine 
space minus an arbitrary union of coordinate hyperplanes, these formulae are 
also the tightest possible upper bounds on the number of isolated roots. We 
also characterize, in terms of sparse resultants, precisely when these upper 
bounds are attained. Finally, we reformulate and extend some of the prior 
combinatorial results of the author on which subsets of coefficients must be 
chosen generically for our formulae to be exact.

Our underlying framework provides a new toric variety setting for 
computational intersection theory in affine space minus an arbitrary union of 
coordinate hyperplanes. We thus show that, at least for root counting, it is 
better to work in a naturally associated toric compactification instead of 
always resorting to products of projective spaces. 
\end{abstract}

\maketitle

\section{Introduction}
\label{sec:intro}
We give a new toric variety context for convex geometric root counts 
for polynomial systems. Our results also improve prior extensions to affine 
space \cite{kho78,dankho,myprize,convexapp,rojaswang,liwang,hsaff} 
of the seminal works \cite{kus75,bernie,kus76,khocompat} on root counting 
in the algebraic torus. In addition to their combinatorial appeal, there has 
been growing excitement about these methods in the computational algebra 
community due to their efficiency and applicability in many industrial problems 
\cite{mobiorobo,emiphd,isawres,dynamiclift}.

\setcounter{footnote}{1}
Let us begin with some notation: Let $E_1,\ldots,E_n$ be nonempty finite
subsets of $\Non$. For any $e\!=\!(e_1,\ldots,e_n)\!\in\!\Non$ let $x^e$
denote the monomial $x^{e_1}_1\cdots x^{e_n}_n$. In this way we will let
$f_1,\ldots,f_n$ be polynomials in the variables $\{x_1,\ldots,x_n\}$ with
(algebraically independent) indeterminate coefficients, such that the set
of exponent vectors occuring in $f_i$ is precisely $E_i$. The set $E_i$ is
called the {\em support}\/ of $f_i$ and this representation specifies
exactly which monomials can appear in $f_i$. All of our root counts will
make maximal use of this monomial term information --- not just the degrees
of the $f_i$. A convenient short-hand will be the following: Let
$E\!:=\!(E_1,\ldots,E_n)$ and $F\!:=\!(f_1,\ldots,f_n)$. Then $E$ is the
{\it support}\/ of $F$ and we call $F$ an {\it $n\!\times\!n$ indeterminate
polynomial system.}\/ We also let $\cC_E$ denote the vector (or sometimes
the set) consisting of all the indeterminate coefficients of all the $f_i$.
If we {\it specialize}\/ some of the coefficients (that is, give them
values chosen from some field)\footnote{We will leave the case of 
more general specializations, e.g., non-trivial polynomial relations amongst 
the coefficients, for future work.} then we say that $F$ has support 
{\em contained}\/ in $E$.

Let $K$ be an algebraically closed field of arbitrary characteristic. For
instance, $K$ can be the complex numbers or the algebraic closure of a
finite field. Also let $|E|$ denote the sum of the cardinalities of the
$E_i$. Our first definition focuses our attention on the generic number of
roots a polynomial system has in a given region $W$, when the
monomial term structure is determined by $E$.

\begin{propdfn}
\label{big}
Let $F$ be an $n\!\times\!n$ indeterminate polynomial system with support
$E$, $\cC_E$ the vector of coefficients of $F$, and $W$ a constructible
subset of $\Kn$. For any $\cC\!\in\!K^{|E|}$ let $\cN_K(E;W;\cC)$ denote
the number of roots of $F|_{\cC_E=\cC}$ lying in $W$, counting
multiplicities.\footnote{See remark \ref{rem:fulton} of section 
\ref{sub:mult} for the definition of intersection multiplicity. } Then there 
exists a {\em proper} algebraic subset
$\Delta\!\subset\!K^{|E|}$, depending on $E$ and $W$, such that
$\cN_K(E;W;\cdot)$ is a constant function on $K^{|E|}\!\setminus\!\Delta$.
We let $\cN_K(E;W)$ denote the value of this constant function. We will
also refer to $\cN_K(E;W)$ as the {\em generic value} of $\cN_K(E;W;\cdot)$
or the {\em generic number of roots} of $F$ in $W$. \qed
\end{propdfn}

\begin{dfn}
Let $[a..b]$ be the set of integers $\{a,a+1,\ldots,b\}$ and for any
(possibly empty) $J\!\subseteq\![1..n]$ define $O_J\!:=\!\{x\!\in\!\Kn\; |
\; x_j\!\neq\!0 \Longleftrightarrow j\!\in\!J\}$. We call $O_J$ an {\em
orbit}. 
\end{dfn}
\noindent
Note that $O_J$ is a relatively open subset of a
$|J|$-dimensional coordinate subspace of $\Kn$. 

Recall that the ($n$-dimensional) mixed volume, $\cM(\cdot)$, takes as
input an $n$-tuple of nonempty compact convex sets in $\Rn$ and always
outputs a nonnegative real number
\cite{bonnie,grunbaum,convexapp,schneider,
polyhomo,isawres,mvcomplex,dynamiclift,drs}. 
\begin{mres}
We will express $\cN_K(E;W)$ in terms of mixed volume for $W$ an arbitrary
union of orbits, $K$ algebraically closed, and any $E$. We will also give a
computational algebraic criterion for precisely when this generic number of
isolated roots is attained, i.e., explicit algebraic equations for $\Delta$. 
Our algebraic criterion is then refined to a more practical computational 
result: a combinatorial classification of the sets of coefficients (subvectors 
of $\cC_E$) whose genericity guarantees that $F$ indeed has exactly 
$\cN_K(E;W)$ isolated roots lying in $W\!$, counting multiplicities.
\end{mres}

The above result is contained in Main Theorems 1--3, the Affine Point
Theorem II, and Corollary \ref{cor:stable}  
of the next section. Examples of our main results
appear in section \ref{sec:examples} and the remaining sections are devoted
to proving our main theorems. Two useful tools applied in our proofs may be
of independent interest: the Antipodality Theorem (\cite{aff} and cf.\
section \ref{sec:mom}) and a
toric variety version of Bernshtein's Theorem (cf.\ section \ref{sub:tv}).
The former tells us how curves behave at toric infinity, while the latter
collects some folkloric facts relating Bernshtein's famous theorem on root
counting \cite{bernie} to intersection theory on toric varieties
\cite{ifulton2,tfulton}.

\section{Summary of Our Main Results}
\label{sec:sum}
We will make the natural restriction of considering only those $E$ for
which $\cN_K(E;W)\!<\!\infty$. Such $E$, which we will call {\em
$W$-nice,}\/ are completely characterized combinatorially in the appendix. It 
will also be helpful to describe certain subspace unions and cones concisely.
\begin{dfn}
\label{dfn:lin}
For any $I\!\subseteq\![1..n]$ let $\hyp(I)\!\subset\!\Kn$ be the union
of coordinate hyperplanes $\bigcup_{j\in I} \{x \; |\\ 
x_j\!=\!0\}$. Also
let $\lin(I)\!\subseteq\!\Rn$ be the coordinate subspace generated by the
subset $\{\hat{e}_j \; | \; j\in I\}$ of the standard basis, and 
let $\sigma_I$ be the cone defined by the intersection of 
$\lin([1..n]\setminus I)$ with 
the nonnegative orthant. Finally, let $\bO$ denote the origin in whatever 
module we work in. In particular, $\hyp(\emptyset)\!=\!\emptyset$ and 
$\lin(\emptyset)\!=\!\sigma_{[1..n]}\!=\!O_\emptyset\!=\bO$.  
\end{dfn}

\subsection{Explicit Formulae}
We give the following recursive formula for $\cN_K(E;W)$. Although perhaps
cumbersome at first glance, our formula contains important intersection
theoretic information that helps extend certain algorithms for solving
polynomial systems \cite{aff} and is also quite practical in low dimensions 
(cf.\  section \ref{sub:mts} and remark \ref{rem:simple}). 
Our result also generalizes, and makes more explicit, an algorithm for 
computing $\cN_\C(E;\Cn)$ (for a smaller class of $E$) alluded to in 
\cite{kho78,dankho}. 
\begin{main}
\label{main:wow}
Let $K$ be any algebraically closed field and suppose
$E\!:=\!(E_1,\ldots,E_n)$ is an $n$-tuple of finite subsets of
$(\N\cup\{0\})^n$ which is nice for $W$, where $W$ is a union of orbits in
$\Kn$. Also, for all
$i,j\in [1..n]$, define $m_{ij}\!:=\!\min\{y_j\; | \; 
(y_1,\ldots,y_n)\!\in\!E_i \}$ and let $m_1,\ldots,m_n$ be the rows of the 
matrix $[m_{ij}]$. Then $\cN_K(E;W)$ is precisely
\[ \sum \limits_{J\subseteq [1..n]} \; 
\sum \limits_{\rho:J^c \hookrightarrow [1..n]} \left[
\left( \prod \limits_{j\in J^c} m_{j\rho(j)} \right)
\cN_K(E_{(J,\rho)};W\cap \lin(\rho(J^c)^c) ) \right], \]
where $(\cdot)^c$ denotes set-theoretic complement within $[1..n]$, 
$E_{(J,\rho)}\!:=\!((E_i-m_i)\cap\lin(\rho(J^c)^c) \; 
| \; i\!\in\!J)$, and $\cN_K(\emptyset;W)$ is defined as $1$ or $0$
according as $W$ is $\bO$ or $\emptyset$. Furthermore, if $W\!=\!\Kni$ for some
$I\!\subseteq\![1..n]$, then $\cN_K(E;W)$ is also the {\em maximum} 
number of isolated roots in $W$, counting multiplicities.
\end{main} 
\begin{rem}
\label{rem:disting}
Our root counting formulae also hold when $E$ is not nice 
for $W,$ provided one counts {\em embedded} \cite[pg.\ 90]{eisenbud} 
zero-dimensional components as well. 
\end{rem}
A simple example of the above formula is given in section \ref{sub:mts} and
its proof appears in section \ref{sub:chow}. Main Theorem 1 is recursive
in the sense that every term on the right-hand side is a  
lower-dimensional or {\em cornered}\/ \cite{rojaswang} case of
$\cN_K(\cdot)$. In particular, the following definition and main result
take care of the ``first'' term $\cN_K(E_{([1..n],\cdot)};W)$. 
\begin{dfn}
Call a $k$-tuple $C\!:=\!(C_1,\ldots,C_k)$ of nonempty subsets of 
$\Rn$ {\em cornered} iff $C_i$ lies in the nonnegative orthant and 
$C_i\cap\{(y_1,\ldots,y_n)\!\in\!\Rn \; | \; y_j\!=\!0\}\!\neq\!\emptyset$ for 
all $i\!\in\![1..k]$ and $j\!\in\![1..n]$. Also, for any 
$a_1,\ldots,a_k\!\in\!\Rn$, define $a\cup C$ to be the $k$-tuple of 
convex hulls $(\conv(\{a_1\}\cup C_1),\ldots,\conv(\{a_k\}\cup C_k))$.  
\end{dfn}
\begin{afpii} 
Fix $I\!\subseteq\![1..n]$ and suppose $E$ is an $n$-tuple of finite
subsets of $(\N\cup\{0\})^n$ which is nice for $\Kni$ and cornered. For
each $i\!\in\![1..n]$ let $a_i\!\in\!E_i\cap\lin(I)$ or set
$a_i\!:=\!\bO$ if $E_i\cap\lin(I)$ is empty. Then
$\cN_K(E;\Kni)=\cM(a\cup E)$ and this generic number is also the {\em
maximum} number of isolated roots in $\Kni$, counting multiplicities. More
generally, if $E$ is instead nice for $O_J$ and cornered, then
$\cN_K(E;O_J)\!=\!\sum_{J'\supseteq J} (-1)^{|J'\setminus J|}
\cM_{J'}$, where $\cM_{J'}$ is the mixed volume corresponding to
the case $W\!=\!\Kn\!\setminus\!\hyp(J')$. 
\end{afpii}
\noindent
The above result is proved in section \ref{sub:embed} and complements 
the author's Affine Point Theorem I which first appeared in \cite{rojaswang}.
\begin{rem}
Note that proposition 1 directly implies that $\cN_K(E;W)$ is additive 
with respect to disjoint unions in $W$. So we can compute 
$\cN_K(E_{([1..n],\cdot)},W)$ 
for general $W$ simply by summing various instances of the Affine Point 
Theorem II. 
\end{rem}
\begin{rem}
Separating into cornered and non-cornered cases also simplifies 
Khovanskii's earlier notion of attached, weakly attached, and strongly 
attached hyperplanes \cite{kho78}.
\end{rem}

The reader need not be alarmed at the prospect of computing an alternating
sum of mixed volumes since a more efficient way to compute $\cN_K(E;O_J)$
is given by the following corollary of Main Theorem 1. This result, which
generalizes a formula for $\cN_\C(E;\Cni)$ due to Huber and Sturmfels
\cite{hsaff}, also seems to yield a more efficient way to compute
$\cN_K(E;W)$ for general $W$ when $n\!>\!2$. 
\begin{cor}
\label{cor:stable}
Following the notation of Main Theorem \ref{main:wow}, fix
$I\!\subseteq\![1..n]$ and suppose $E$ is nice for $\Kni$. Then
$\cN_K(E;\Kni)=\cSM_{I^c}(E)$. Furthermore, if
$J\!\subseteq\![1..n]$, $E$ is instead $O_J$-nice, and $\Omega$ is a stable
subdivision of $E$, then $\cN_K(E;O_J)\!=\!\sum_C \cM(C)$, where the sum is 
over all stable cells $C\!\in\!\Omega$ such that the inner normal of the 
lifted cell $\widehat{C}$ has support $J^c$.
\end{cor}
\noindent
The quantity $\cSM_J(E)$ is a new convex geometric entity called the {\em
$J$-stable mixed volume}. We refer the reader to \cite{hsaff} for its
definition, and to \cite{polyhomo,hsaff} for the definitions of
subdivisions, lifted cells, and stable cells. The support of a vector is
simply the set of indices corresponding to its nonzero coordinates.
Corollary \ref{cor:stable} is proved in section \ref{sub:chow}.

Better still, we can determine precisely when our formulae count 
the number of roots exactly, even when some of the coefficients are fixed and 
only a few coefficients are generic. 

\subsection{Algebraic and Combinatorial Criteria for Exactness} For any
$w\!\in\!\Rn$, let $E^w_i$ be the set of points $y\!\in\!E_i$ which minimize
the standard inner product $w\cdot y$. Similarly, for any polytope
$P\!\subset\!\Rn$, let $P^w$ denote the face of $P$ with inner normal
$w$. Also let $E^w\!:=\!(E^w_1,\ldots,E^w_n)$ and recall that a {\em facet}\/ 
is a polytope face of codimension 1. A key innovation of Bernshtein's seminal
work on root counting is the algebraic condition he gave for his formula  
to be the {\em exact}\/ number of roots. Sadly, this ``second half'' of 
Bernshtein's Theorem is not sufficiently explored in the literature. So we 
give the following generalization, proved in section \ref{sub:algcond}.
\begin{main} 
\label{main:res}
Following the notation of definitions \ref{big} and \ref{dfn:lin} and Main 
Theorem \ref{main:wow}, suppose $W\!=\!\Kni$ and that the coefficients of $F$ 
have all been specialized to constants in $K$. Let $S$ be the polytope 
$\sum^n_{i=1}\conv(E_i)$. Then the following condition implies that $F$ has 
exactly $\cN_K(E;\Kni)$ roots, counting multiplicities, in $\Kni$:
\begin{itemize}
\item[(a$_2$)]{$\prod \res_{E^w}(F)\!\neq\!0$, where the product is over all 
unit inner facet normals $w\!\in\!\Rni$ of $S$, and }
\item[(b$_2$)]{if $n\!>\!1$ then for all $J\!\subsetneqq\![1..n]$ 
containing $I$, and all injections $\rho : J^c \hookrightarrow [1..n]$ such 
that $\rho(J^c)\cap I\!=\!\emptyset$ and $\prod_{j\in J^c} m_{j\rho(j)}>0$, 
\[ \cN_K\left(E_{(J,\rho)};\lin(\rho(J^c)^c)\cap(\Kni);\cC_E\right) = 
\cN_K\left(E_{(J,\rho)};\lin(\rho(J^c)^c)\cap(\Kni)\right).\]}
\end{itemize}
Furthermore, the converse implication holds as well if 
$\cN_K(E_{([1..n],\cdot)};\Kni)\!>\!0$. In particular, (a$_2$) and (b$_2$) 
together imply that the zero set of $F$ in $\Kni$ is zero-dimensional or empty. 
\end{main}
\noindent 
The sparse resultant $\res_*(F)$ is described at length in 
\cite{gkz90,chowprod,jandi1,
combiresult,gkz94,isawres} and our notation is explained in section 
\ref{sub:algcond}. Sharper criteria for the cases 
$\cN_K(E_{([1..n],\cdot)};\Kni)\!=\!0$ are also discussed in section 
\ref{sub:algcond}. 
\begin{rem}
Parallel to Main Theorem 1, Main Theorem 2 is also recursive since for any 
$\vartheta\!\subsetneqq\![1..n]$, $\lin(\vartheta)\cap(\Kni)$ can be naturally 
identified with the complement of a union of coordinate hyperplanes in 
$K^{|\vartheta|}$. 
\end{rem}

Alternatively, we can give combinatorial criteria for exactness 
which are always sufficient and necessary. More precisely, let $c_{i,e}$ 
denote the (indeterminate) coefficient of the $x^e$ term of $f_i$. If an 
$n$-tuple 
$D\!:=\!(D_1,\ldots,D_n)$ satisfies $D_i\!\subseteq\! E_i$ for all 
$i\!\in\![1..n]$ then we simply abbreviate this as $D\!\subseteq\!E$. For any 
such $D$ define $\cC_D\!=\!\{c_{i,e} \; | \; i\!\in\![1..n]$, $e\!\in\!D_i\}$.
\begin{dfn}
\label{dfn:count} 
We say that $D$ {\em $W$-counts} $E$ iff (0) $D\!\subseteq\!E$, (1) $D$ 
and $E$ are nice for $W$, and (2) for any specialization over $K$ of the
coefficients $\cC_E\!\setminus\!\cC_D$, a generic specialization of the
remaining coefficients $\cC_D$ suffices to make $F$ have exactly
$\cN_K(E;W)$ roots lying in $W$, counting multiplicities.
\end{dfn}
\noindent
So by proposition \ref{big} we at least know that $E$ always 
$W$-counts $E$ if $E$ is $W$-nice. 

Define $\supp(D)\!:=\!\{i \; | \; D_i\!\neq\!\emptyset\}$,
$D\cap\lin(J)\!=\!(D_1\cap\lin(J),\ldots,D_n\cap\lin(J))$, and $D\cap
E^w\!:=\!(D_1\cap E^w_1,\ldots,D_n\cap E^w_n)$. Our final main theorem
gives an exact convex geometric criterion for when $D$ $(\Kni)$-counts $E$. 
\begin{main} 
\label{main:count}
Following the notation of definitions \ref{big} and \ref{dfn:lin} and Main 
Theorem \ref{main:wow}, suppose that $W\!=\!\Kni$ and let 
$S$ be the polytope $\sum^n_{i=1}\conv(E_i)$. Then $D$ $(\Kni)$-counts $E 
\Longleftrightarrow$ one of the following exclusive conditions holds: 
\begin{enumerate}
\item{$\cN_K(E;\Kni)\!=\!0$ and for all $J\!\subseteq\![1..n]$ 
containing $I$, $\supp(D\cap\lin(J))$ contains a subset essential for 
$E\cap\lin(J)$.} 
\item{$\cN_K(E;\Kni)\!>\!0$ and 
\begin{itemize}
\item[(a$_3$)]{for each face of $S$ with an inner normal $w\!\in\!\Rni$, 
pick a single such $w$. Then 
for all of these $w$, $\supp(D\cap E^w)$ contains a subset essential for 
$E^w$, and }
\item[(b$_3$)]{if $n\!>\!1$ then for all $J\!\subsetneqq\![1..n]$ 
containing $I$, and all injections $\rho : J^c \hookrightarrow [1..n]$ such 
that $\rho(J^c)\cap I\!=\!\emptyset$ and $\prod_{j\in J^c} m_{j\rho(j)}>0$, 
the $|J|$-tuple  
\[ ((D_i-m_i)\cap \lin(\rho(J^c)^c) \; | \; i\!\in\!J) \]  
$W_{(J,\rho)}$-counts $E_{(J,\rho)}$, where 
$W_{(J,\rho)}\!:=\!\lin(\rho(J^c)^c)\cap(\Kni)$.}
\end{itemize}} 
\end{enumerate}
\end{main}
The definition of essentiality, which is a combinatorial geometric
condition, appears in the appendix. We thus obtain a recursive combinatorial 
condition for when the zero set of $F$ in $\Kni$ consists of exactly  
$\cN_K(E;\Kni)$ points, counting
multiplicities. Our final main theorem is proved in section \ref{sub:algcond} 
as well. Here we deal mainly with genericity conditions for {\em global}\/ 
root counting, so we will leave the classification 
of $O_\vartheta$-counting (when $\vartheta\!\neq\![1..n]$) for another paper. 

\section{Examples}
\label{sec:examples}
In the following examples, any mixed volume computation will follow easily
(even by hand) from the definition or basic properties of the mixed volume
\cite{bonnie,grunbaum,schneider}. In particular, it useful to recall
the following formula for the $n\!=\!2$ case: 
$\cM(P_1,P_2)\!=\!\area(P_1+P_2)-\area(P_1)-\area(P_2)$.
\subsection{Comparisons to the Generalized B\'{e}zout Theorems}
\label{sub:bez}
Although mixed volume bounds can be hard to compute for some extremely
large polynomial systems, they do have the advantage that they are always
at least as good as any B\'ezout-type bound. Also, current mixed volume
software is already fast enough to have been useful in many industrial
problems, e.g., \cite{emiphd,dynamiclift}. Here we will give an example of
a family of polynomial systems whose mixed volume root counts are
significantly better than any generalized B\'ezout bound.

However, let us first recall what is meant by a generalized B\'ezout bound.
A good reference is \cite{wam92} so we will only quickly outline the most
general (zero-dimensional) version of B\'ezout's Theorem: Given a partition
of $\{x_1,\ldots,x_n\}$ into sets of cardinality $n_1,\ldots,n_\lambda$, the 
corresponding {\em multihomogeneous}\/ B\'ezout Theorem gives an explicit 
formula for $\cN_K(E;\Pro^{n_1}_K\!\times \cdots \times\!\Pro^{n_\lambda}_K)$ 
as a polynomial expression involving the degrees of the $f_i$ with 
respect to the chosen sets of variables.\footnote{Polynomial roots in a 
toric compactification are 
described in sections \ref{sec:mom} and \ref{sub:embed}, and are formalized 
in definition \ref{dfn:shift} and lemma \ref{lemma:shift}. In particular, 
products of projective spaces are special cases of toric compacta. } 
Implicit in the grouping of variables chosen is an embedding $\Kn 
\hookrightarrow \Pro^{n_1}_K\!\times \cdots \times\!\Pro^{n_\lambda}_K$ and in 
this way we obtain an upper bound on $\cN_K(E;\Kn)$.

One can then try to group variables so that this method gives as tight
an upper bound on $\cN_K(E;\Kn)$ as possible, but the following example
shows that this bound can be very loose, no matter how one groups
variables.
\begin{ex}{\bf (Spiky Newton Polytopes)} Consider the indeterminate
polynomial system $F\!:=\!(c_{10}+c_{11}x^d_1+ c_{12}x_2+\cdots
+c_{1n}x_n,\ldots,c_{n0}+c_{n1}x^d_1+c_{n2}x_2+\cdots + c_{nn}x_n)$, where
$d\!\in\!\N$ and $n\geq 2$. Clearly, the Newton polytopes of $F$ are all
identical and equal to the ``spike''
$P\!:=\!\conv(\bO,d\hat{e}_1,\hat{e}_2,\ldots,\hat{e}_n)$. The Affine Point
Theorem II then tells us that $\cN_K(E;\Kn)=\cM(P,\ldots,P)=n!\V(P)=d$,
where $E$ is the support of $F$.

However, the usual B\'ezout Theorem tells us that $\cN_K(E;\Kn)\!\leq\!d^n$.
Can this be significantly improved by going to a multi-homogeneous version? 
The answer is ``yes, but not enough:'' the best one can do is
$\cN_K(E;\Kn)\!\leq\!nd$. This bound can be obtained by using two groups of
variables: $\{x_1\}$ and $\{x_2,\ldots,x_n\}$. That this is the best one
can do with any multihomogeneous version of B\'ezout's Theorem is most
easily proved geometrically: It is easy to see that computing the optimal
generalized B\'ezout bound is equivalent to finding a product of scaled
standard simplices, with smallest volume, which contains $P$. (This
reduction is described more explicitly, but in a different context, in
\cite{myaverage}.) 
\end{ex}

More generally, one can use Main Theorem 3 to determine when a particular
B\'ezout Theorem generically matches (or exceeds) a mixed volume root
count: One simply lets $E$ be the $n$-tuple of vertex sets of the
corresponding products of simplices, and lets $D$ be the $n$-tuple of
vertex sets of the Newton polytopes in question. {}From there, one checks
the corresponding {\em counting}\/ criterion (cf.\ definition 
\ref{dfn:count}). 

\subsection{Generic Local Intersection Multiplicity} 
\label{sub:mult}
Setting $W\!=\!\bO$ in
Main Theorem 1 we immediately obtain a method for computing the {\em
generic}\/ intersection multiplicity, at the origin, of a general sparse
system of $n$ polynomials in $n$ unknowns. An alternative general
algorithm, potentially more efficient in higher dimensions, is the special
case $J\!=\!\emptyset$ of Corollary \ref{cor:stable}. For example, if $F$
is a $2\!\times\!2$ polynomial system with {\em cornered}\/ support $E$, we
obtain from the Affine Point Theorem II that
$\mu(\bO;F)\!=\!\cN_K(E;\bO)\!=\!\cM(\bO\cup E)-\cM(E)$ for generic
$\cC_E$. This last formula already generalizes an earlier result of  
Warren \cite[Theorem 3]{warren} for the {\em unmixed}\/ case ($E_1\!=\!E_2$) 
over $\C$. 

However, it is important to note that $\cN_K(E;\bO)$ is {\em not}\/, in
general, an upper bound on intersection multiplicity at the origin: For 
example, it is 
easily verified that the polynomial system $(x+y^2,x+x^2+y^2)$ has an
isolated root at $\bO$ with multiplicity $4$. (One simply notes that this
system has no roots other than $\bO$ and concludes by B\'ezout's Theorem in
$\Pro^2_K$.) However, setting $E\!:=\!(\{(1,0),(0,2)\},\{(1,0),(2,0),
(0,2)\})$, our last paragraph implies that 
$\cN_K(E;\bO)\!=\!\cM(\bO\cup E)-\cM(E)\!=\!4-2\!=2$.
\begin{rem}
\label{rem:fulton}
More generally, if $F$ is an $n\!\times\!n$ polynomial system over $K$, 
then the intersection multiplicity $\mu(\zeta,F)$ of a zero-dimensional 
component $\zeta\!\in\!\Kn$ of $Z(F)$ can be defined as the dimension 
of the quotient ring $\{\frac{r(x)}{s(x)}\!\in\!K(x) \; | \; \gcd(r,s)\!=\!1, 
\; s(\zeta)\!\neq\!0\}/\langle F \rangle$ as a $K$-vector space 
\cite[Example 7.1.10 (b)]{ifulton2}. For the purposes of definition 
\ref{big}, we will set $\mu(\zeta,F)\!=\!+\infty$ when $\zeta$ lies 
in a {\em positive}-dimensional component of $Z(F)$. 
\end{rem}

\subsection{Our Main Theorems in Two Dimensions} 
\label{sub:mts}
\begin{figure}[bth]
\begin{picture}(500,140)(-45,-10)

\put(0,0){
\begin{picture}(133,100)(0,0)
\put(18,29){$E_1$ ($\supseteq\!\supp(f_1)$)} 
\put(-10,0){\line(1,0){100}} \put(80,-7){$e_1$}
\put(0,-10){\line(0,1){100}} \put(-10,84){$e_2$}
\put(0,20){\circle*{3}} \put(2,19){$\alpha_1$}
 \put(0,20){\line(0,1){20}}
\put(0,40){\circle*{3}} \put(-12,40){$\alpha_2$}
 \put(0,40){\line(2,3){20}}
\put(10,50){\circle{3}} \put(5,44){$\alpha_3$}
\put(20,50){\circle{3}} \put(22,48){$\alpha_4$}
\put(20,70){\circle*{3}} \put(8,70){$\alpha_5$}
 \put(20,70){\line(2,1){20}}
\put(40,80){\circle{3}} \put(41,81){$\alpha_6$}
 \put(40,80){\line(-2,-3){40}}
\end{picture}}

\put(108,0){
\begin{picture}(133,100)(0,0)
\put(3,45){$E_2$ ($\supseteq\!\supp(f_2)$)} 
\put(-10,0){\line(1,0){100}} \put(80,-7){$e_1$}
\put(0,-10){\line(0,1){100}} \put(-10,84){$e_2$}
\put(10,0){\circle{3}} \put(8,-9){$\beta_1$}
 \put(10,0){\line(0,1){10}} 
\put(20,0){\circle*{3}} \put(18,-9){$\beta_2$}
 \put(20,0){\line(2,1){40}}
\put(10,10){\circle*{3}} \put(5,14){$\beta_3$}
 \put(10,10){\line(2,1){80}}
\put(20,10){\circle{3}} \put(22,9){$\beta_4$}
\put(60,20){\circle{3}} \put(61,16){$\beta_5$}
 \put(60,20){\line(1,1){30}}
\put(60,30){\circle{3}} \put(51,24){$\beta_6$}
\put(70,30){\circle{3}} \put(72,27){$\beta_7$}
\put(90,50){\circle*{3}} \put(92,51){$\beta_8$}
\end{picture}}

\put(216,0){
\begin{picture}(133,100)(0,0)
\put(10,7){$S\!=\!\conv(E_1+E_2)$} 
\put(-10,0){\line(1,0){100}} \put(80,-7){$e_1$}
\put(0,-10){\line(0,1){100}} \put(-10,84){$e_2$}
\put(10,20){\circle{3}}
 \put(10,20){\line(1,0){10}}
\put(20,20){\circle{3}}
 \put(20,20){\line(2,1){40}}
\put(60,40){\circle{3}}
 \put(60,40){\line(1,1){30}}
\put(90,70){\circle{3}}
 \put(90,70){\line(2,3){40}}
\put(130,130){\circle{3}}
 \put(130,130){\line(-2,-1){100}}
\put(30,80){\circle{3}}
 \put(30,80){\line(-2,-3){20}} 
\put(10,50){\circle{3}}
 \put(10,50){\line(0,-1){30}} 
\put(110,120){\circle{3}} \put(110,110){\circle{3}} \put(110,100){\circle{3}}
\put(100,110){\circle{3}} \put(100,100){\circle{3}} \put(90,100){\circle{3}}
\put(90,90){\circle{3}} \put(80,100){\circle{3}} \put(90,80){\circle{3}}
\put(80,90){\circle{3}} \put(90,70){\circle{3}} \put(80,80){\circle{3}}
\put(80,70){\circle{3}} \put(70,80){\circle{3}} \put(60,90){\circle{3}}
\put(70,70){\circle{3}} \put(60,80){\circle{3}} \put(50,90){\circle{3}}
\put(60,70){\circle{3}} \put(50,80){\circle{3}} \put(70,50){\circle{3}}
\put(60,60){\circle{3}} \put(40,80){\circle{3}} \put(60,50){\circle{3}}
\put(40,70){\circle{3}} \put(30,80){\circle{3}} \put(60,40){\circle{3}}
\put(40,60){\circle{3}} \put(30,70){\circle{3}} \put(40,50){\circle{3}}
\put(30,60){\circle{3}} \put(30,50){\circle{3}} \put(20,60){\circle{3}}
\put(20,50){\circle{3}} \put(20,40){\circle{3}} \put(10,50){\circle{3}}
\put(20,30){\circle{3}} \put(10,40){\circle{3}} \put(20,20){\circle{3}}
\put(10,30){\circle{3}} \put(10,20){\circle{3}}
\dashline[50]{2}(10,30)(90,70)
\dashline[50]{2}(10,30)(50,90)
\end{picture}}

\end{picture} 

\caption{The sets $E_1$, $E_2$, and $E_1+E_2$, along with their underlying 
convex hulls. (The points of $E_1$ and $E_2$ are labelled according to their 
corresponding polynomial coefficients.) The subdivision of the Minkowski sum 
shows that, in this example, $\cM(E)$ can actually be expressed as a single 
determinant. } 
\label{fig:useful}

\end{figure}
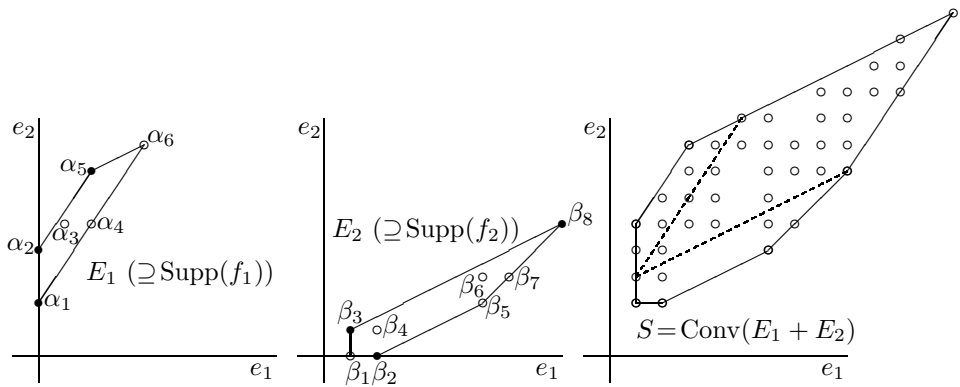

Let $n\!=\!2$ and consider the following bivariate polynomial system: 
\begin{eqnarray*}
f_1(x,y) & := & \alpha_1 y^2 + \alpha_2 y^4 + \alpha_3 xy^5 + \alpha_4 x^2y^5 +
\alpha_5 x^2y^7 + \alpha_6 x^4y^8 \\
f_2(x,y) & := & \beta_1 x + \beta_2 x^2 + \beta_3 xy + \beta_4 x^2y +
\beta_5 x^6y^2 + \beta_6 x^6y^3 + \beta_7 x^7y^3 + \beta_8 x^9y^5 
\end{eqnarray*}
How do we get a tight upper bound on the number of isolated affine roots
of $F\!:=\!(f_1,f_2)$? One way is to set $E\!:=\!\supp(F)$ and apply the Affine
Point Theorem I \cite{rojaswang}. In which case we obtain that $F$ has no
more than $\cM(\bO\cup E)\!=\!53$ isolated roots. (It is also easily
verified that the best generalized B\'ezout bound is $(\deg_x f_1)(\deg_y
f_2)+(\deg_x f_2)(\deg_y f_1)\!=\!92$.) However, it is clear that $E$ is
nice for $K^2$ (by lemma \ref{lemma:wow} of the appendix) and not cornered, so 
let us see if Main Theorem 1 can do better. Figure \ref{fig:useful} below 
clarifies the preceding (and upcoming) calculations. 

Following the notation of Main Theorem 1, we obtain 
$m_1\!=\!(0,2)$, and $m_2\!=\!(1,0)$. So Main Theorem 1 asserts that
\begin{eqnarray*}
\cN_K(E;K^2) & \!=\! & \cN_K(E_{\{1,2\}};K^2)+
\cN_K(\{\bO,(0,2)\},\{0\}\!\times\!K) + 
2\cN_K(\{\bO,(1,0)\},K\!\times\!\{0\}) + 2\cN_K(\emptyset,\bO)\\ 
& \!=\! & 32 + \cN_K(\{0,2\},K) + 2 \cN_K(\{0,1\},K) + 2. 
\end{eqnarray*}
(The last equality follows from the Affine Point Theorem II.) The remaining
unknown terms add up to $4$ (by the Affine Point Theorem II again, or simply
the fundamental theorem of algebra), so we finally obtain the tight upper 
bound $\cN_K(E;K^2)\!=\!38$. 

We can also give a precise algebraic condition for when this $F$ has {\em
exactly}\/ thirty eight affine roots, counting multiplicities. Applying
Main Theorem 2 (and figure \ref{fig:useful}) to our example,   
we see that the only $w$ we need worry about in condition (a$_2$) are (in 
counter-clockwise order) $(-1,2)$, $(-1,-1)$, $(-3,2)$, $(1,-2)$, and 
$(3,-2)$. Furthermore, the corresponding sparse resultants are easily seen to 
be $1$, $1$, $1$, $\alpha_5\beta_8-\alpha_6\beta_3$, and $1$. (The 
papers \cite{morea,chowprod,combiresult} and the book \cite{gkz94} 
contain some very nice examples of how to compute low-dimensional sparse
resultants.) Condition (b$_2$) then clearly specializes to two one-dimensional 
cases of condition (a$_2$). (More conservatively, the fundamental theorem of 
algebra could also be applied to (b$_2$).) So it is not hard to see that 
condition (b$_2$) is equivalent to $\alpha_2$ and $\beta_2$ being nonzero. 
Since each individual term of the summation from Main Theorem 1 was positive, 
we thus obtain that $F$ has exactly thirty eight affine roots, counting 
multiplicities, {\bf iff} 
$(\alpha_5\beta_8-\alpha_6\beta_3)\alpha_1\alpha_2\beta_2\beta_3
\beta_8\!\neq\!0$.

\begin{figure}[bth]
\begin{picture}(410,70)(-35,-5) 

\put(-20,-50){
\begin{picture}(133,100)(0,0)
\put(41,80){\circle*{3}} \put(38,85){$C_1$} \put(82,80){\circle*{3}}
\put(79,85){$C_2$} \put(45,50){$\{1\},\{2\}$}
\end{picture}}

\put(108,-50){
\begin{picture}(133,100)(0,0)
\put(26,70){\line(1,1){30}}
\put(26,70){\circle*{3}} \put(56,100){\circle*{3}} \put(26,85){$C_1$}
\put(72,70){\line(1,1){30}}
\put(72,70){\circle*{3}} \put(102,100){\circle*{3}} \put(72,85){$C_2$}
\put(50,50){$\{1,2\}$}
\end{picture}}

\put(256,-50){
\put(26,100){\line(1,-1){30}}
\put(26,100){\circle*{3}} \put(56,70){\circle*{3}} \put(43,85){$C_1$}
\put(72,70){\line(0,1){30}}
\put(72,70){\circle*{3}} \put(72,100){\circle*{3}} \put(74,85){$C_2$}
\begin{picture}(133,100)(0,0)
\put(45,50){None}
\end{picture}}

\end{picture}

\caption{The essential subsets for 3 different pairs of plane polygons. 
(The segments in the middle pair are meant to be parallel.) } 
\label{fig:myfirst}

\end{figure}
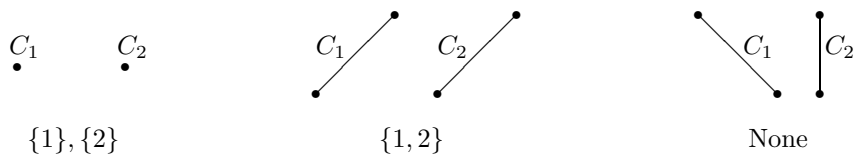

Similarly, by Main Theorem 3 (and making use of figures 
\ref{fig:useful} and \ref{fig:myfirst} (cf.\ the appendix)) we obtain that the 
genericity of $\cC_D$ implies $F$ has exactly thirty eight affine roots 
(counting multiplicities) $\Longleftrightarrow \cC_D$ contains at least one 
coefficient from {\em each}\/ of the following sets: 
$\{\alpha_1,\beta_2\}$, $\{\alpha_1\}$,
$\{\alpha_1,\beta_5\}$, $\{\alpha_1\}$, $\{\alpha_8\}$,
$\{\alpha_6,\beta_8\}$, $\{\alpha_5,\alpha_6\}$,$\{\beta_3,\beta_8\}$,
$\{\alpha_5,\beta_3\}$, $\{\beta_3\}$, and $\{\alpha_2,\beta_3\}$ (in 
counter-clockwise order, from condition (a$_3$)); and $\{\alpha_2\}$ and
$\{\beta_2\}$ (from condition (b$_3$)). For example, 
regardless of how the other eight coefficients have been
specialized, it suffices to choose the vector 
$(\alpha_1,\alpha_2,\alpha_5,\beta_2,\beta_3, \beta_8)\!\in\!K^6$ 
generically for $F$ to have exactly thirty eight roots (counting
multiplicities). In other words, $D$ $K^2$-counts $E$, where 
$D=(\{(0,2),(0,4),(2,7)\},\{(2,0),(1,1),(9,5)\})$. The dark points 
in figure \ref{fig:useful} represent $D$.

\section{Background and Terminology}
Aside from a few variations, we will follow the same notation as
\cite{convexapp,isawres,rojaswang,hsaff,dynamiclift}. In those papers one
can also find some of the definitions below described at a more leisurely
pace. We will also liberally quote, e.g., from
\cite{grunbaum,algeom,sha,clo,schneider}, various simple facts from
convex and algebraic geometry that we will use. However, for the
convenience of the reader, we will review a few notions.

For any $q_1,\ldots,q_n\!\in\!\Rn$, let $[q_1,\ldots,q_n]$ denote the
$n\!\times\!n$ matrix whose $i^\thth$ column is $q_i$. It will be useful to
recall the following facts concerning the Smith normal form of an
integral matrix \cite[Chap.\ 3.7]{jacobson1} 
\begin{propdfn}
\label{prop:sln}
\cite{unimod1,unimod2}
An {\em integral basis}\/ for $\Rn$ is a vector space basis for $\Rn$ which
is also a $\Z$-module basis for $\Zn$. Equivalently, a basis
$\{u_1,\ldots,u_n\}$ for $\Rn$ is integral iff
$[u_1,\ldots,u_n]\!\in\!\glnz$. Such a basis {\em respects}\/ a rational
subspace $L$ of $\Rn$ iff $\{u_1,\ldots,u_{\dim L}\}$ is a $\Z$-module
basis for $L\!\cap\!\Zn$. Furthermore, given any rational bases for $L$ and
a complementary subspace, an integral basis for $\Rn$ respecting $L$ can be
found within $O(n^{11})$ arithmetic operations. \qed 
\end{propdfn}
Defining $x^\cU\!:=\!(x^{u_{11}}_1\cdots
x^{u_{n1}}_n,\ldots,x^{u_{1n}}_1\cdots x^{u_{nn}}_n)$ for any
$n\!\times\!n$ matrix $\cU\!:=\![u_{ij}]$ with integer entries, we then see
that the above proposition tells us how we can find a {\em monomial}\/ change 
of variables which converts a polynomial into a form involving as few
variables as possible.

We will use $\supp(f)$ and $\newt(f)$ for, respectively, the support and
Newton polytope (the convex hull of $\supp(f)$) of any $f\!\in\!K'[x^{\pm
1}_1,\ldots,x^{\pm 1}_n]$, where $K'\!:=\!K[\Lambda]$ and $\Lambda$ is any set 
of algebraically independent indeterminates. 
\begin{dfn}
For any {\em weight} $w\!\in\!\Rn$ the {\em initial term
polynomial} $\init_w(f)(x)$ is $\sum_{e \in \supp(f)^w} c_e x^e$. More
generally, if $B\!\supseteq\!\supp(f)$, we define the {\em relativized}
initial term polynomial $\init_{w,B}(f)(x)\!:=\! \sum_{e \in B^w \cap
\supp(f)} c_e x^e$. Also, any $c_e$ with $e$ on the boundary of $\newt(f)$
is called a {\em boundary coefficient} of $f$. 
\end{dfn}
\noindent
Alternatively, when $w\!\in\!\Zn$, we can simply substitute $x \mapsto t^w
x\!:=\!(t^{w_1}x_1,\ldots,t^{w_n}x_n)$ into $f$ and define $\init_w(f)$ as
the coefficient of the term of {\em lowest}\/ degree in $t$.

More generally, for any $k\times n$ polynomial system $F$ (with constant
or indeterminate coefficients), we define the {\em initial term system}\/
$\init_w(F)$ to be $(\init_w(f_1),\ldots,\init_w(f_k))$. Also, if a 
$k$-tuple $C\!:=\!(C_1,\ldots,C_k)$ satisfies $C_i \! \supseteq \!
\supp(f_i)$ for all $i\!\in\![1..k]$, then we say that $C$ {\em contains the 
support of $F$}\/ and we define the {\em relativized}\/ initial term system
$\init_{w,C}(F)$ to be $(\init_{w,C_1}(f_1),\ldots,\init_{w,C_k}(f_k))$. An
especially important property of initial term systems is the following.
\begin{prop}
\cite{convexapp}
\label{prop:over}
Suppose $F$ is an $n\!\times\!n$ indeterminate polynomial system with support 
$E\!=\!(E_1,\ldots,E_n)$. In particular, we assume that each $E_i$ is nonempty. 
Then for generic $\cC_E$ and any $w\!\neq\!\bO$, the polynomial system 
$\init_{w,C}(F)$ has {\em no} roots in $\Ksn$. \qed 
\end{prop}
\noindent
Note that for any polynomial system $F$ with support contained in $C$, the
set $\{\init_{w,C}(F) \; | \;
w \!\in\! \Sn \}$ is finite: When $\supp(f_i)\!=\!C_i$ for all $i$, we can 
construct a bijection between the set of initial term systems and the face 
lattice of $\conv(\sum C_i)$, simply by picking a single inner normal $w$ for 
each face.

There is a rich interplay between
the combinatorial geometric structure of $\newt(f)$ and the topology of the
zero set of $f$ and we will see again (in section \ref{sec:mom} and beyond)
that initial term polynomials are extremely valuable in this respect.

\subsection{Algebraic Geometry}
\label{sub:algebraic}
As usual, we will let $Z(F)$ denote the zero scheme of $F$ in $\Kn$. We
will make some use of algebraic cycles (e.g., finite formal
$\Z$-linear combinations of closed subvarieties of some toric variety),
rational equivalence, and intersection theory, so let us also recall the
following facts and definitions \cite{algeom,sha,ifulton1,ifulton2,tfulton}: 

\begin{enumerate}
\item{For any cycle $\cA$, $\supp(\cA)$ is the union of all
closed subvarieties $V$ such that the coefficient of $V$ within $\cA$ is
nonzero. Also, a divisor is said to be {\em effective}\/ iff all its
coefficients are nonnegative. }
\item{There is a natural intersection product ``$\cap$'' on the group
of all cycles on a variety $\cX$ giving this group a (commutative) ring
structure called the {\em Chow ring}\/ of $\cX$, $\chow(\cX)$. This product
is also compatible with {\em rational equivalence} 
\cite[pp.\ 10, 15--17]{ifulton2}.}
\item{A {\em $0$-cycle}\/ on $\cX$ is a cycle of the form
$\cD\!=\!\sum n_\zeta \{\zeta\}$ where each $\zeta$ is a point and 
$n_\zeta\!\in\!\Z$. When $\cX$ is complete, the homorphism from the group of 
$0$-cycles on $\cX$ to 
$\Z$ defined by $\sum n_\zeta \{\zeta\}\mapsto\sum n_\zeta$ is invariant under 
rational equivalence and is called the {\em degree map}, $\deg(\cdot)$. } 

\item{Any intersection, $\cD$, of $\dim \cX$ many divisors in $\chow(\cX)$ 
is rationally equivalent to some $0$-cycle, and thus has a well-defined ({\em
cycle class}\/) degree. Furthermore, if each divisor is effective, the
coefficient of any zero-dimensional component $\zeta$ of such an
intersection is its {\em intersection multiplicity}, or {\em intersection
number}, $\mu(\zeta;\cD)$ \cite[Example 7.1.10 (b)]{ifulton2}.} 

\end{enumerate}

The most advanced prerequisites we will require from algebraic geometry
will be the belief in certain theorems dealing with divisor intersections
on toric varieties. Good general references are \cite{dannie,oda,ifulton2,
tfulton}. The toric variety facts we'll need are covered in the next
section so we now state the main intersection theoretic result we'll use.

\begin{thm}
\label{thm:intersect} \cite[Theorem 12.2]{ifulton2} Suppose $\cT$ is a
complete $n$-dimensional variety over an algebraically closed field, and
$\cD_1,\ldots,\cD_n$ are effective Cartier divisors such that each line
bundle $\cO(\cD_i)$ is generated by its sections. Also let $\cD$ denote the
intersection product of these divisors in $\chow(\cT)$. Then the
intersection number of any distinguished component of $\cD$ is nonnegative.
Furthermore, if we let $\cI$ denote the sum of the intersection numbers of
the distinguished components of $\cD$, then \[\cI\!\leq\!\deg \cD,\] and
equality holds if $\cD$ is zero-dimensional or empty. \qed
\end{thm}
\begin{rem}
See \cite{ifulton2} for the definition of a distinguished component. For
our purposes, suffice it to say that a zero-dimensional irreducible component 
is a distinguished component, and the converse holds as well when $\cD$ is
zero-dimensional or empty. Also, embedded \cite[pg.\ 90]{eisenbud} 
zero-dimensional components are distinguished components. 
\end{rem}
\noindent
Precise conditions for equality in the above inequality are subtle and
difficult to find in the literature. However, we conjecture that equality
always holds in the cases where we will apply this theorem. This has
already been verified in a particular case, giving a refinement of
B\'ezout's Theorem over $\C$ \cite{mikebezout}.

\subsection{Toric Varieties}
\label{sub:toric}
We will assume the reader to be familiar with fans and the
construction of toric varieties from fans and polytopes. Excellent
references are \cite{kkms,dannie,oda,tfulton,gkz94,grobandpoly}.

Let $T\!:=\!\Ksn$, which is sometimes called the {\em algebraic torus}.

\begin{dfn}
\label{dfn:polyfan}
Let $P\!\subset\!\Rn$ be an $n$-dimensional rational polytope. We will
associate to $P$ its {\em (inner) normal fan} $\fan(P)$ as follows: The rays of
this fan are generated by the inner facet normals of $P$, and to each (not
necessarily proper) face $P^w$ of $P$ we associate the cone $\sigma_w$
generated by the rays corresponding to the facets containing $P^w$. Each
$\sigma_w$ is also called a {\em (inner) normal cone} of $P$. 
\end{dfn}
\noindent
It is useful to think of the {\em duals} of the cones of $\fan(P)$ as
``angle'' cones. In fact, it easy to show that for any $w$ there is a small
ball $B\!\subset\!\Rn$, centered at the origin, such that
$B\cap\sigma^\vee_w\!=\!B\cap (P-v)$, for some $v\!\in\!\relint P^w$.

We will be working with the following class of toric varieties. 
\begin{dfn}
Following the notation of definition \ref{dfn:polyfan}, we will let $\cT_P$
be the toric variety over $K$ corresponding to the normal fan of $P$. We
call $\cT_P$ the {\em toric compactification} of $T$ corresponding to $P$.
\end{dfn}
It follows that $\cT_P$ is $n$-dimensional, rational, projective, normal,
integral, separated, and complete \cite{tfulton}. The toric variety $\cT_P$
also has a naturally embedded copy of $T$ (cf.\ theorem \ref{thm:orbit}). For 
certain $P$ the toric variety $\cT_P$ is also nonsingular but we will not need 
this fact. 
We will also say that any point of $\cT_P\!\setminus\!\Ksn$ is {\em at 
infinity} and sometimes refer to $\cT_P\!\setminus\!\Ksn$ as {\em toric 
infinity}. Since our polynomial systems will have a priori specified supports, 
$F$ will usually have far fewer extraneous roots in an appropriately chosen 
$\cT_P$ than in $\Pro^n_K$. Hence toric compactifications are the spaces where 
we will actually be counting roots of polynomial systems. Toward this end, it
will be useful to recall the correspondence between the topology of $\cT_P$  
and the face structure of $P$. However, we will need a little more 
notation before stating this correspondence as a theorem.
\begin{dfn}
\label{dfn:toric}
Given any $w\!\in\!\Rn$, we will use the following notation: 
\begin{itemize}
\item[$U_w=$]{$\spec(K[x^e \; | \; e\!\in\!\sigma^\vee_w\cap\Zn])=$The affine 
chart of $\cT_P$ corresponding to the cone
$\sigma_w$ of $\fan(P)$} 
\item[$L_w=$]{The $\dim(P^w)$-dimensional subspace of $\Rn$ parallel to the
face $P^w$ of $P$} 
\item[$x_w=$]{The point in $U_w$ corresponding to the semigroup
homomorphism $\sigma^\vee_w\cap\Zn \longrightarrow \{0,1\}$ mapping $p
\mapsto \delta_{w\cdot p,0}$, where $\delta_{ij}$ denotes the Kronecker
delta} 
\item[$O_w=$]{The $T$-orbit of $x_w=$ The $T$-orbit corresponding to $\relint 
P^w$} 
\item[$V_w=$]{The closure of $O_w$ in $\cT_P$} 
\item[$p_w=$]{The first integral point {\em not} equal to $\bO$ met along the 
ray generated by $w$ (when $w\!\in\!\Qns$) } 
\end{itemize}
\end{dfn}
\noindent
Note that $L_w$ is a face of the cone $\sigma_w^\vee$ so $x_w$ is indeed
well-defined. Also, recalling that a closed point in an affine toric
variety can be identified with a semigroup homomorphism \cite[Chap.\
1.3]{tfulton}, it is clear that any point $x\!\in\!O_w$ is completely
determined by $w$ and the (nonzero) values of $x(\cdot)$ on any $\Z$-module
basis of $L_w\cap \Zn$. Note that our characterization of $x_w$ is a
slight variation of that of \cite{tfulton} but is easily seen to be
equivalent. In particular, our $x_w$ is the same as Fulton's $x_\sigma$
when $\sigma\!=\!\sigma_w$.
\begin{ex} {\bf (Certain Cornered Polytopes)} Suppose $P$ is an $n$-dimensional
rational polytope with a vertex ${\mathrm v}$ such that the edges
emanating from ${\mathrm v}$ generate the nonnegative orthant as a cone.
Suppose further that the coordinates of $w$ are all {\em non}\/negative.
Then $O_w\!\cong\!O_J$, where $J\!=\!\supp(w)^c$. In
particular, we see that $\dim P^w=\dim O_w=n-|\supp(w)|$ and $\dim \sigma_w
= |\supp(w)|$. Note that here $x_w$ is the 0-1 vector with support
$\supp(w)^c$. Furthermore, $U_{(1,\ldots,1)}\!\cong\!\Kn$ 
so we can thus conclude that $\Kn$ embeds naturally within such a $\cT_P$.
This example will be especially important in our approach to affine root
counting. 
\end{ex}
\begin{ex}
\label{ex:basic}
Suppose $w,w'\!\in\!\Rn$. Then, relative to $\cT_P$, the defining ideal
$I_{w'}\!\subset\!K[x^e \; | \; e \in \sigma^\vee_w \cap \Zn]$ of $V_{w'}
\cap U_w \subset U_w$ is $K[x^e \; | \; e\in
\!(\sigma^\vee_w\!\setminus\!L_{w'})\!\cap\!\Zn]$. 
\end{ex}

With our orbit notation in place, we can now state the following important
result.
\begin{thm}
\label{thm:orbit}
\cite[Chap.\ 3.1]{tfulton} 
The toric variety $\cT_P$ is the disjoint union $\coprod O_w$, where a
single inner normal $w$ is chosen for each (not necessarily proper) face of
$P$. Also, for all $w\in \Rn$,
\begin{enumerate}
\item{$\dim O_w = \dim V_w = \dim P^w$}
\item{$U_w=\coprod O_v$, where a single inner normal $v$ is chosen for each
(not necessarily proper) face containing $P^w$. In particular, $U_w$ is
always an $n$-dimensional open subvariety of $\cT_P$. }
\item{$O_w\!=\!\spec(K[x^e \; | \;
e\!\in\!L_w\!\cap\!\Zn])\!\cong\!(\Ks)^{\dim P^w}$.} 
\item{If $d\!=\!\dim
P^w$ then $V_w$ is isomorphic to any toric compactification of $(\Ks)^d$
corresponding to a polytope $Q\subset\R^d \subseteq\Rn$ which is
$\glnz$-similar to a translate of $P^w$. \qed}
\end{enumerate}
\end{thm}
\noindent 
In particular, there is an order-preserving correspondence 
between the faces (resp.\ face interiors) of $P$ and the 
orbit closures (resp.\ orbits) of $\cT_P$. Also, there is an 
order-{\em reversing}\/ correspondence between the affine 
charts of $\cT_P$ and the faces of $P$. The above result 
is also contained in \cite{ksz92,gkz94} but in the setting 
where $\cT_P$ is defined via an explicit projective embedding.  

Since toric compactifications will be the spaces in which we analyze the
roots of $F$, it will be useful to embed the support of $F$ within an
$n$-tuple of nonempty integral polytopes $\cP\!:=\!(P_1,\ldots,P_n)$ and 
define $P$ as a function of $\cP$.
We can then consider the roots of $F$ within $\cT_P$ as follows: Each
(nonzero) polynomial $f_i$ defines a Weil divisor $\divisor(f_i)$ in
$\cT_P$ \cite[Chap.\ 3.3]{tfulton}. The closure (in $\cT_P$) of the zero
scheme of $f_i$ in $\Ksn$ is a summand of $\divisor(f_i)$ and is the
portion of $\divisor(f_i)$ we are actually interested in. To isolate this
portion of $\divisor(f_i)$ we will add another specially defined divisor
(depending on $P$ and $P_i$) to $\divisor(f_i)$. This will cancel out the
negative part of $\divisor(f_i)$ but sometimes introduce extraneous
components. In any case, the zero scheme of $F$ in $\Ksn$ is thus embedded
in an intersection of effective divisors in $\cT_P$. In section
\ref{sub:chow} we will show
how to eliminate some of these extra components, modify $P$ so that  
$\cT_P$ has a naturally embedded copy of $\Kni$, and thus derive our method
for affine root counting. The construction of our divisors is detailed in
the next section. 

\section{The Importance of Roots at Toric Infinity} 
\label{sec:mom}
Here we point out two, more or less folkloric results on toric divisors.
Combined with the Antipodality Theorem \cite{aff}, these two results 
considerably simplify the proof of our toric compactification version 
(theorem \ref{thm:toric} in the next subsection) of the BKK bound 
\cite{bernie}. 

First we give the following definition to help us find the right $\cT_P$,
and the right divisor to add to $\divisor(f_i)$, for our root counting
theory to go through.
\begin{dfn}
\label{dfn:impo1}
Let $Q\!\subset\!\Rn$ be an {\em integral} polytope. We will say that a fan
$\cF$ is {\em compatible}\/ with $Q$ iff every normal cone of $Q$ is a
union of cones of $\cF$. We will also say that a {\em rational} polytope
$P\!\subset\!\Rn$ is {\em compatible} with $Q$ iff $\fan(P)$ is compatible 
with $Q$. Also, following the notation of definition \ref{dfn:toric}, we 
define the integer $\gamma_w(Q)\!:=\!-\min
\limits_{v\in Q} \{ v\!\cdot\! p_w\}$ for any $w\!\in\!\Qns$. 
\end{dfn}
\begin{ex}
\label{ex:compat}
It is easily shown that $\sum P_i$ is always compatible with $P_1,\ldots,P_n$. 
Compatibility was applied earlier in \cite{khocompat,tfulton} and the
terminology ``sufficiently fine decomposition'' was used in the first
reference. 
\end{ex}

Next we describe precisely which divisors we will be intersecting. 
\begin{dfn}
\label{dfn:impo2}
Assuming $P\!\subset\!\Rn$ is a rational polytope compatible with 
an integral polytope $Q\!\subset\!\Rn$, let $\cE_P(Q)\!:=\!\sum 
\gamma_{w}(Q) V_w$, where $w$
ranges over all the inner facet normals of $P$. We call $\cE_P(Q)$ {\em the
torus-invariant divisor of $\cT_P$ corresponding to $Q$.} Also, 
set $\cD_P(0,Q)\!:=\!\cT_P$ and, for any polynomial $f$ with 
$\emptyset\!\neq\!\supp(f)\!\subseteq\!Q$, define 
$\cD_P(f,Q)\!:=\!\divisor(f)+\cE_P(Q)$. 
\end{dfn}
\noindent
It follows by definition that $\cD_P(f,Q)$ is always effective
\cite{tfulton} and invariant under translations of $P$ and identical
translations of $Q$ and $f$. This turns out to be good for root counting in
$\Ksn$ but bad for root counting in $\Kn$. Hence we will modify the
definition of $\cD_P(f,Q)$ in section \ref{sub:chow}. 
\begin{dfn}
\label{dfn:intersect} 
Suppose $F\!=\!(f_1,\ldots,f_k)$ is a $k\!\times\!n$ polynomial system
over $K$ with support contained in a $k$-tuple of nonempty integral polytopes
$\cP\!=\!(P_1,\ldots,P_k)$. Then $P\!\subset\!\Rn$ is {\em compatible} 
with $\cP \Longleftrightarrow P$ is compatible with $P_1,\ldots,P_k$. 
Furthermore, when this is the case, we define $\cD_P(F,\cP)$ to be the 
intersection product $\bigcap^k_{i=1} \cD_P(f_i,P_i)\!\in\!\chow(\cT_P)$. 
\end{dfn}
It is easy to see that even as schemes 
$\Ksn\cap Z(F)\!=\!\Ksn\cap\cD_P(F,\cP)$ when $\dim P\!=\!n$. 
Also note that $\cD_P(f,Q)$ is precisely the closure $\overline{\Ksn\cap
Z(f)}$ in $\cT_P$ if $\newt(f)\!=\!Q$ and $\dim P\!=\!n$. 
By our last observation, we could just set $P_i\!:=\!\newt(f_i)$ for 
all $i$ in the construction of $\cD_P(F,\cP)$ in order 
to work directly with $\overline{\Ksn\cap Z(F)}$. However, this is not always
advantageous computationally and it actually behooves us to fully
understand the cases where $\newt(f_i)$ is not compatible with $P$ or
$\newt(f_i)\!\neq\!P_i$. One reason is that for precise {\em sparse}\/
affine root counting, it is necessary to know precisely what happens to 
$\divisor(f_i)$ as lots of coefficients of $F$ are specialized to $0$.

So let us now find explicitly the behavior of $\cD_P(F,\cP)$ within a
neighborhood of toric infinity. The following lemma, which is a direct
consequence of the development followed in \cite{tfulton} or \cite{gkz94}, 
shows us that $U_w\cap\cD_P(F,\cP)$ can be described by a relatively simple 
ideal.
\begin{lemma}
\label{lemma:once}
Let $w\!\in\!\Rn$. Then, following the notation of definitions \ref{dfn:toric} 
and \ref{dfn:intersect}, the defining ideal in 
$K[x^e \; | \; e\in\sigma^\vee_w\cap \Zn]$ of $U_w\cap \cD_P(F,\cP)$ is 
$\langle x^{b_1}f_1,\ldots,x^{b_k}f_k \rangle$, for any 
$b_1,\ldots,b_k\!\in\!\Zn$ such that $b_i+P^w_i\subseteq L_w$ for all 
$i\!\in\![1..k]$. 
\qed 
\end{lemma}
\noindent
Should one be so inclined, the intersection multiplicity of a component of
$\cD_P(F,\cP)$ can be computed by restricting to an appropriate chart $U_w$ 
and this lemma gives one an explicit coordinate ring to work in.

The following is a more computational version of the above lemma and is
easily proved by localization.
\begin{cor}
\label{cor:cool}
Following the notation of lemma \ref{lemma:once}, the underlying
topological spaces of
$O_w\cap\cD_P(F,\cP)$ and
\[\spec(K[x^e \; | \; e\in L_w\cap \Zn] / \langle \init_{w,b_1+P_1}(x^{b_1}
f_1),\ldots,\init_{w,b_k+P_k}(x^{b_k} f_k) \rangle)\]
are homeomorphic.

In particular, if $\{u_1,\ldots,u_n\}\!\subseteq\!L_w\cap \Zn$ is a
generating set, then the map defined by $x \mapsto (x(u_1),\ldots,x(u_n))$
is an isomorphism from $O_w$ onto a subvariety of $\Ksn$, and
$[z\!\in\!O_w\cap\cD_P(F,\cP) \Longleftrightarrow \init_{w,\cP}(F)$ 
vanishes at the point $(z(u_1),\ldots,z(u_n))]$. \qed 
\end{cor}

It is a frequent misconception that the last equivalence should begin with
$z\!\in\!O_w\cap\overline{\Ksn\cap Z(F)}$. This is false unless, for
instance, $P_i\!=\!\newt(f_i)$ for all $i$. (A simple counterexample is
$(n,P,f,w)\!:=\!(1,[0,2],x^2+x,1)$.)
Also, it is extremely important to note that the intersection multiplicity
of a component of $\cD_P(F,\cP)$ lying in $O_w$ can {\em not}\/ always be 
determined by this corollary when $w\!\neq\!\bO$. (Simply observe the 
case $P_1\!=\!P_2\!=\!P\!=\!3\conv\{\bO,\hat{e}_1,\hat{e}_2\}$,
$F\!=\!(x^3_1-x^2_1,x^3_2-x^2_2)$, and $w\!:=\!(1,1)$.) Thus there is a loss 
of information as we intersect with $O_w$ and pass from $F$ to its initial 
term systems.

The last two results thus tell us that the relativized initial term systems
$\init_{w,\cP}(F)$ (for $w\!\neq\!\bO$) describe the topological behavior
of a particular divisor intersection (canonically defined by $F$, $\cP$, 
and $P$) at a piece of toric infinity. This generalizes the classical 
construction of how the terms of highest total degree depict the closure of 
the zero scheme of $F$ at the projective hyperplane at infinity.

Our toric variety $\cT_P$ also gives us an interesting way to detect 
excess components in the zero set of $F$ in $\Ksn$. 
\begin{anti} \cite{aff}
Suppose $P\!\subset\!\Rn$ is an $n$-dimensional rational polytope, $Y$ is a
curve in $\Ksn$, and $\overline{Y}$ is the closure of $Y$ in $\cT_P$. Also,
identify $\{O_w\!\subset\!\cT_P \; | \; w\!\neq\!\bO\}$ with a partition of
$\Sn$ into {\em cells} by letting $w$ and $w'$ in $\Sn$ belong to the same
cell iff $O_w\!=\!O_{w'}$. Then
\begin{enumerate}
\item{$\overline{Y}$ must touch two (possibly identical) cells whose union
does not lie in any closed coordinate hemisphere. } 
\item{If $\overline{Y}$ touches a cell intersecting some open hemisphere $\cH$
then it must also touch a cell intersecting $\Sn\!\setminus\!\cH$. }
\end{enumerate}
Furthermore, $\overline{Y}$ must touch at least two topologically separated
cells. \qed
\end{anti}
\begin{rem}
The cells described above are in fact stereographic images of the facets of 
the dual (or polar) of $P$. 
\end{rem}

\subsection{A Toric Variety Version of Bernshtein's Theorem} \label{sub:tv}
Before stating our generalization of Bernshtein's Theorem, we point out a
very useful immediate corollary of the Antipodality Theorem and corollary 
\ref{cor:cool}. 
\begin{cor}
\label{cor:crux}
Following the notation of corollary \ref{cor:cool}, $\cD_P(F,\cP)$ has positive
dimension $\Longrightarrow \init_{w,\cP}(F)$ has a root in $\Ksn$ for some
$w\!\neq\!\bO$. \qed
\end{cor}
\noindent
D. N. Bernshtein proved the case of corollary \ref{cor:crux} where
$\Ksn\cap\cD_P(F,\cP)$ is positive-dimensional and $K\!=\!\C$ \cite{bernie}. 
His proof used an ingenious Puiseux series construction, but unfortunately
Puiseux series expansions are not always defined for algebraic curves over
a field of {\em positive}\/ characteristic. Hence our need for antipodality
theorems.

By combining the following lemma
with theorem \ref{thm:intersect}, we see that working within $\cT_P$ allows
us to reduce the computation of intersection numbers (generically) to the
evaluation of a mixed volume.

\begin{lemma}
\label{lemma:cartier} \cite[Chap.\ 5.4]{tfulton} Following the notation of
definitions \ref{dfn:impo2} and \ref{dfn:intersect}, the cycle class degree of 
$\cE_P(P_1) \cap \cdots \cap \cE_P(P_n) \in \chow(\cT_P)$ is precisely 
$\cM(E)$. Furthermore, for all $i$, the line bundle $\cO(\cD_P(f_i,P_i))$ is 
generated by its
sections. \qed 
\end{lemma}

Putting all our machinery together, we can derive the following toric
variety version of the BKK bound.
\begin{thm}
\label{thm:toric}
Suppose $F$ is an $n\!\times\!n$ polynomial system over $K$ with support
contained in an $n$-tuple $\cP$ of nonempty integral 
polytopes in $\Rn$. Further suppose that $P\!\subset\!\Rn$ is an 
$n$-dimensional rational polytope compatible with $\cP$. Then the 
zero scheme of $F$ in $\Ksn$ embeds naturally as a subscheme of the toric 
cycle $\cD_P(F,\cP)$ and:  
\begin{enumerate}
\item{If $\cD_P(F,\cP)$ is zero-dimensional or empty then $\cD(F,\cP)$ 
consists of exactly $\cM(\cP)$ points, counting multiplicities.} 
\item{If $\cD_P(F,\cP)$ is
positive-dimensional and $\cM(\cP)\!=\!0$ then $\cD_P(F,\cP)$ has {\em no}
zero-dimensional irreducible components in $\cT_P$.} 
\item{If $\cD_P(F,\cP)$ is positive-dimensional and $\cM(\cP)\!>\!0$ then 
$\cD_P(F,\cP)$ has strictly less than
$\cM(\cP)$ zero-dimensional irreducible components in $\cT_P$, counting
multiplicities.}
\end{enumerate}
\end{thm}
\begin{rem}
Assertion (3) appears to be new for the case $\ch K\!\neq\!0$.
The case $(K,P_1,\ldots,P_n)\!=\!(\C,\newt(f_1),\ldots, \newt(f_n))$ first
appeared in \cite{bernie} and was stated as a root count over $\Csn$
instead of $\cT_P$. Assertions (1) and (2) (over a general algebraically 
closed field) then appeared implicitly in \cite{dannie}, but this was not as 
well known as it should have been. 
\end{rem}
\noindent
{\bf Proof:} That $\Ksn\cap Z(F)$ embeds as a subscheme of $\cD_P(F,\cP)$ 
follows immediately from our previous observations regarding definitions 
\ref{dfn:impo2} and \ref{dfn:intersect}. 

Assertion (1) then follows 
immediately from theorem \ref{thm:intersect} and lemma \ref{lemma:cartier} 
since $\cD_P(f_i,P_i)$ and $\cE_P(P_i)$ are rationally equivalent. Assertion 
(2) follows similarly since the intersection multiplicity of a zero-dimensional 
irreducible component of $\cD_P(F,\cP)$ is positive \cite{ifulton2}. Fulton 
stated this concise argument in \cite{tfulton} for the 
case $K\!=\!\C$ and his proof has the added benefit that it is independent of 
the (algebraically closed) field where one is working.

For the last assertion (3) we will generalize a novel homotopy proof due to
D. N. Bernshtein \cite{bernie}. First note that if any $f_i$ is identically
zero then there can be no zero-dimensional components and we are done. Thus
we may assume that no $f_i$ is identically zero. Our generalization of
Bernshtein's argument can then be outlined as follows: 
\begin{itemize}
\item[(i)]{Pick a point $y\!\in\!O_w$, for some $w\!\in\!\Rns$, which lies in
a positive-dimensional component of $\cD_P(F,\cP)$.} 
\item[(ii)]{Construct a generic polynomial system $G$ with $n$-tuple of
Newton polytopes $\cP$ and distinguished root $z\!\in\!\Ksn$ such that
$F(z)\!\neq\!\bO$.}
\item[(iii)]{Construct a rational algebraic curve
$L\!\subset\!\cT_P\!\times\!\Pro^1_K$
with parameterization $\overline{l} : \Pro^1_K \longleftrightarrow L$ such
that $\overline{l}(0)\!=\!y$ and $\overline{l}(1)\!=\!z$.} 
\item[(iv)]{For all $i\!\in\![1..n]$, define
$h_i(x,t)\!\in\!K[t,x^{\pm 1}_1,\ldots,x^{\pm 1}_n]$ to be the polynomial
obtained by clearing denominators from the reduced form of the rational
function
$f_i(x)g_i(\overline{l}(t))-g_i(x)f_i(\overline{l}(t))$. Show that
$h_i(x,0)$ (resp.\ $h_i(x,1)$) is a nonzero scalar multiple of $f_i(x)$
(resp.\ $g_i(x)$). }
\item[(v)]{Let $H(x,t)\!:=\!(h_1(x,t),\ldots,h_n(x,t))$ and consider the
subscheme $Z\!:=\!\overline{(\Ksn\!\times\!K)\cap Z(H)}$ of
$\cT_P\!\times\!\Pro^1_K$. Show that $\cD_P(F,\cP)\cong
(\cT_P\!\times\!\{0\})\cap Z$. }
\item[(vi)]{Show that the natural $(n+1)^\st$ coordinate projection defined
on $\Ksn\!\times\!K$ extends to a proper morphism $\pi : 
\cT_P\!\times\!\Pro^1_K \longrightarrow \Pro^1_K$ with
$\pi(L)\!=\!\Pro^1_K$. }
\item[(vii)]{Define $\cY$ to be the union of all $1$-dimensional irreducible
components of $Z$ with surjective image under $\pi$. Show that the support
of the zero-dimensional part of $\cD_P(F,\cP)$ is contained in
$\supp((\cT_P\!\times\!\{0\})\cap \cY)$. } 
\item[(viii)]{Show that $\cY\!\cap\!(\cT_P\!\times\!\{0\})$ consists of exactly
$\cM(\cP)$ points, counting multiplicities. } 
\end{itemize}
Assuming the above steps, (3) then follows immediately since $y\!\in\!L$,
$L\!\subseteq\!\cY$, and thus
$y\!\in\!\supp(\cY\!\cap\!(\cT_P\!\times\!\{0\}))$, i.e., the
zero-dimensional part of $\cD_P(F,\cP)$ consists of strictly fewer than 
$\cM(\cP)$ points, counting multiplicities.

To complete our proof, we now proceed to prove each individual step.

\noindent
{\bf (i):} Easy, by the Antipodality Theorem.

\noindent
{\bf (ii):} By generic we will
specifically mean that $f_1,\ldots,f_n$ are all nonzero at all roots of $G$
in $\Ksn$ and that $G$ has exactly $\cM(\cP)$ roots in $\Ksn$ (counting
multiplicities). That such $G$ occur generically follows easily from
proposition \ref{prop:over}, corollary \ref{cor:crux}, assertion (1) (which
we've already proved), and the fact that the intersection of any two
generic conditions is again a generic condition.

\noindent
{\bf (iii):} We will first construct the parameterization $\overline{l}$
and then the corresponding complete curve $L$.

Since $y$ is completely determined by the (nonzero) values of $y(\cdot)$ on
any $\Z$-module basis of $L_w\cap\Zn$, let $\{u_1,\ldots,u_n\}$ be any
basis for $\Zn$ respecting $L_w\!\cap\!\Zn$. (Such a basis is guaranteed to
exist by proposition \ref{prop:sln}.) Let $\cU\!:=\![u_1,\ldots,u_n]$ and
$[v_1,\ldots,v_n]\!:=\!\cV\!:=\!\cU^{-1}$. By assumption $\cU\!\in\!\glnz$
so clearly $\cV\!\in\!\glnz$. Now let $l : \Ks \longrightarrow \Ksn$ be the
following parameterization of a toric line: \[ ((1-t)y^\cU + t z^\cU)^\cV
\]
where, quite naturally, $y^\cU\!:=\!(y(u_1),\ldots,y(u_n))$ and $t$ is a
new variable. Then it is easily verified that $l(1)\!=\!z$ (via the general
identity $(x^A)^B\!=\!x^{AB}$). Note that $l$ naturally defines a rational
function from $\Pro^1_K$ to $\cT_P\!\times\!\Pro^1_K$ via $[t\!:\!1]
\mapsto (l_1(t),\ldots,l_n(t))\!\times\![t\!:\!1]$. So by \cite[Prop.\
2.1]{silverman} this rational function extends uniquely to a morphism
$\overline{l}$. This is our desired $\overline{l}$ and, of course,
$\overline{l}(1)\!=\!z$.

Now let $L$ be the closure in $\cT_P\!\times\!\Pro^1_K$ of the subvariety of 
$(\Ks)^{n+1}$ defined by the ideal
\[ \langle (1-t)y(u_1) + tz^{u_1} - x^{u_1}, \ldots, (1-t)y(u_n) + tz^{u_n}
- x^{u_n} \rangle. \] 
It is then easily verified
that the hyperplane $\cT_P\!\times\!\{t_0\}$ intersects $L$ in the unique
point $l(t_0)\!\in\!\Ksn$ for all but finitely many $t_0\!\in\!K$. (Simply
solve the resulting binomial equation by exponentiating by $\cV$.) It is
also clear that $\cT_P\!\times\!\{t_0\}$ does not meet $L$ within
$\Ksn\!\times\!\{t_0\}$ for the
remaining values of $t_0$. Thus $L$ must indeed be a curve. Since $L$ is
closed, it must also be complete and equal to the graph of $\overline{l}$
in $\cT_P\!\times\!\Pro^1_K$. By corollary \ref{cor:cool} it easily follows
that $L\!\cap\!O_{(w,1)}\!=\!(z,0)$ (for $w$ as in (i)) and thus
$\overline{l}(0)\!=\!z$.

\noindent
{\bf (iv):} First note that the definition makes sense since we can just
substitute $l$ for $\overline{l}$. Now to verify that $h_i(x,t)$ satsifies
our desired properties, note that $l(t)^{u_i}\!=\!(1-t)y(u_i)+tz^{u_i}$ by
a straightforward calculation. If we define $d\!:=\!\dim L_w$ and write
$e\!=\!\alpha_1 u_1 + \cdots + \alpha_n u_n$, it then becomes clear that
$\ord_t(l(t)^e)\!=\!\alpha_{d+1}+\cdots+\alpha_n$. By changing the signs of
the columns of $\cU$ where necessary, we can then assume that
$w\!\cdot\!u_i\!\geq\!0$ for all $i\!\in\![1..n]$ and write
\[f_i(l(t))\!=\!t^{b_i}(\init_{w,P_i}(f_i)|_{x=(y(\pi_w(\hat{e}_1)),\ldots,
y(\pi_w(\hat{e}_n)))}+ \mathrm{higher \ order \ terms \ in \ } t)\]
\[g_i(l(t))\!=\!t^{b_i}(\init_{w,P_i}(g_i)|_{x=(y(\pi_w(\hat{e}_1)),\ldots,
y(\pi_w(\hat{e}_n)))}+ \mathrm{higher \ order \ terms \ in \ } t)\] where
$b_i\!\in\!\Z$ and $\pi_w : \Rn \longrightarrow L_w$ is the natural
projection defined by the basis $\cU$. In particular, note that
$t^{-b_i}f_i(l(t))|_{t=0}\!=\!0$ (by corollary \ref{cor:cool} and the
definition of $l$) and $t^{-b_i}g_i(l(t))|_{t=0}\!\neq\!0$ (by corollary
\ref{cor:cool} and the definition of $G$). Thus $h_i(x,t)\!=\!t^{-b_i} \kappa(t)
(f_i(x)g_i(l(t))-g_i(x)f_i(l(t)))$, for some $\kappa(t)\!\in\!K(t)$ satisfying
$\ord_t(\kappa)\!=\!0$. So we are done.

\noindent
{\bf (v):} First note that for any $\omega\!\in\!\R^{n+1}$, there is a
$b\!\in\!\Z^{n+1}$ such that $\newt(x^b h_i)$ touches every facet of the
cone $\sigma^\vee_\omega$. (This follows easily since all the $f_i$ are
nonzero and $G$ has $n$-tuple of Newton polytopes $\cP$.) Thus by
definition \ref{dfn:impo2},
$\overline{(\Ksn\!\times\!K)\cap Z(h_i)}\!=\!\cD_{P\times
[0,1]}(h_i,P_i\times [0,d_i])$, where $d_i$ is the $t$-degree of $h_i$.
Lemma \ref{lemma:once} then implies that
for any $w'\!\in\!\Rns$, the defining ideal of $Z\cap U_{(w',1)}$ is
generated by $x^{b_1}h_1,\ldots,x^{b_n}h_n$, for suitable
$b_1,\ldots,b_n\in\Zn$. By theorem \ref{thm:orbit} we know that
$\cT_P\!\times\!\{0\}\cong V_{\hat{e}_{n+1}}$, so by lemma \ref{lemma:once}
the defining ideal of
$(\cT_P\!\times\!\{0\})\cap U_{(w',1)}$ is principal and generated by $t$.
Since $h_i\!\equiv\!f_i \pmod{t}$, we thus see that the defining ideal of
$Z\!\cap\!(\cT_P\!\times\!\{0\})\!\cap\!U_{(w',1)}$ is generated by $t$ and
$x^{b_1}f_1,\ldots,x^{b_n}f_n$. Lemma \ref{lemma:once} then implies that
$\cD_P(F,\cP)\!\cap\!U_{w'}$ has a defining ideal with generators
$x^{b_1}f_1,\ldots,x^{b_n}f_n$. Patching together charts, we are done.

\noindent
{\bf (vi):} This follows easily from \cite[Chap.\ 2.4]{tfulton}.

\noindent
{\bf (vii:)} Clearly, any zero-dimensional irreducible component $\zeta$ of
$\cD_P(F,\cP)$ must be contained in some positive-dimensional irreducible 
component
$\Sigma$ of $Z$. If $\dim \Sigma\!>\!1$ then $\zeta$ must lie in a
positive-dimensional component of $\Sigma \cap (\cT_P\!\times\!\{0\})$ ---
a contradiction. Thus $\dim \Sigma\!=\!1$. Clearly $0\!\in\!\pi(\Sigma)$
and $\pi(\Sigma)\!\neq\!\{0\}$
so by properness we must have $\pi(\Sigma)\!=\!\Pro^1_K$. So $\Sigma$ must
be a component of $\cY$.

\noindent
{\bf (viii:)} Let $t_0\!\in\!\Pro^1_K$. If $\cY$ is irreducible then it
follows directly from the definition of intersection multiplicity that
$\sum \mu(\zeta)$, where $\zeta$ ranges over all zero-dimensional
irreducible components of $\cY\cap(\cT_P\times\{t_0\})$, is precisely the sum
of the ramification indices \cite{silverman} of $\pi^{-1}(t_0)$. The latter
number in turn is precisely the degree of the map $\pi$ and is thus
independent of the point $t_0$ \cite[Examples 4.3.7 and
7.1.15]{ifulton2}.\footnote{Note that in our definitions of ramification
index and degree, we are including the inseparability degree. This is
relevant when the characteristic of $K$ is positive.} If $\cY$ is reducible
then we can extend the definition of degree simply by summing the degrees
of $\pi|_{\cY_j}$ over all the irreducible components $\cY_j$ of $\cY$ and then
our preceding identity still holds.

Thus it suffices to compute $\sum \mu(\zeta)$ for {\em any}\/ $t_0$. In
particular, by construction,
we already know that this number is precisely $\cM(\cP)$ when $t_0\!=\!1$. \qed

The above theorem is quite useful for root counting in $\Ksn$ but still has
the nagging problem that it doesn't give the exact number of roots when the 
intersections are ill-behaved --- more precisely, when $\cD_P(F,\cP)$ 
intersects toric infinity.
However, our theorem (when combined with corollary \ref{cor:cool}) at least
provides us with a computational method for knowing exactly when this
happens. (Indeed, Main Theorem 2 is based on this very fact!) Also, when 
$\cD_P(F,\cP)$ is zero-dimensional, the precise number of 
roots, counting multiplicities, can {\em still}\/ be obtained as follows.
\begin{cor}
\label{cor:subtract}
Following the notation of theorem \ref{thm:toric}, suppose further that
$\cD_P(F,\cP)$ is zero-dimensional or empty. Let $\cI\!:=\!\sum \mu(\zeta)$ 
where the sum ranges over all components $\zeta$ of
$\cD_P(F,\cP)\!\setminus\!\Ksn$. Then the number of roots of 
$F$ in $\Ksn$ is precisely $\cM(\cP)\!-\!\cI$, counting multiplicities. \qed
\end{cor}
\noindent
This approach to {\em exact}\/ (as opposed to generic) root counting is
pursued further in \cite{lamina,resvan} and was independently 
suggested in \cite{msw95} (in the special case of multihomogeneous systems)
and \cite[pp.\ 180--185 and 215--216]{janphd} (not counting some intersection 
multiplicities). 

Intuitively, it is a weaker condition to require $\cD_P(F,\cP)$ to be
zero-dimensional than to require all the roots of $F$ to be isolated {\em
and}\/ lie in $\Ksn$. This statement is made more precise in the next
section and in section \ref{sub:algcond} we will also give a combinatorial
characterization of the stronger hypothesis.

Another natural question which still remains is how to extend our analysis
to other spaces --- for example, $\Kn$. We do this in the next section.
\begin{rem}
We point out that {\bf all} of the results of this section hold for more
general complete toric varieties as well --- in particular, toric varieties
corresponding to (compatible) complete fans. The modifications are minor
and we have omitted them simply because toric compacta corresponding to
polytopes are sufficiently powerful for our particular root counting
problems in affine space. 
\end{rem}
\begin{rem}
\label{poly}
A more elementary (but longer) proof of parts (1) and (2) of theorem
\ref{thm:toric} can be obtained by generalizing Huber and Sturmfels' proof
of Bernshtein's Theorem \cite{polyhomo} to arbitrary algebraically closed
fields: One first replaces the use of theorem \ref{thm:intersect} and lemma
\ref{lemma:cartier} in our intersection theoretic proof by the
combinatorics of {\em mixed subdivisions}. Then, the use of Puiseux series
in their proof is replaced by some algebraic curve theory \`a la the proof of 
(viii). The resulting polyhedral proof requires no more machinery than that 
already used in the proof of part (3). However, for the sake of brevity we 
will omit this alternative proof.
\end{rem}

\section{Proofs of Our Five Main Results} \label{sec:proofs}
We now expand our applications of toric compacta to root counting in affine
space. We will begin by proving the Affine Point Theorem II and then
proceed to prove Main Theorem 1, Corollary \ref{cor:stable}, and Main
Theorems 3 and 2. 

\subsection{Affine Embeddings}
\label{sub:embed}
Contrary to what one might expect, a toric compactification $\cT_P$ does
{\em not}\/ always contain a naturally embedded copy of $\Kn$. This
technicality forces us to require $P$ to satisfy an additional hypothesis
before we apply $\cT_P$ to root counting in $\Kn$. The following definition
is the first of our two main tricks for applying toric intersection theory
to affine root counting.

\begin{dfn}
\label{dfn:compat}
We say a rational polytope $P\!\subset\!\Rn$ is {\em cornered} iff
$\fan(P)$ contains the nonnegative orthant as one of its cones. More
generally, for any $I\!\subseteq\![1..n]$, $P$ is {\em $I$-cornered} iff 
$\sigma_I$ is one of the cones of $\fan(P)$ (following the 
notation of definition \ref{dfn:lin}). 
\end{dfn}
\noindent 
Note that cornering is different for polytopes and $k$-tuples of 
point sets: For polytopes, $\emptyset$-cornering is easily seen to be 
equivalent to a translate of $P$ being identical to the nonnegative orthant in 
a neighborhood of $\bO$. For a $k$-tuple $(C_1,\ldots,C_k)$, cornering refers 
to the {\em position}\/ of each $C_i$ within the nonnegative orthant 
$\sigma_\emptyset$.

Our last definition is well-motivated for the following reason. 
\begin{prop}
\label{prop:corn}
If $P\!\subset\!\Rn$ is $I$-cornered then $\cT_P$ has a naturally embedded
copy of $\Kni$. More precisely, for such a $\cT_P$, $\Kni\!\cong\!U_w$,
where $w$ is the 0-1 vector with support $I^c$. \qed
\end{prop}

We now show how to construct a special $I$-cornered $\Pc$ from any given 
$k$-tuple of polytopes in $\Rn$. 
\begin{alg}
\label{alg:corner} \ 
\begin{itemize}
\item[Input:]{A positive integer $n$, a $k$-tuple of {\em nonempty} integral 
polytopes $\cP\!=\!(P_1,\ldots,P_k)$ lying in the nonnegative orthant of $\Rn$, 
and a subset $I\!\subseteq\![1..n]$.} 
\item[Output:]{An $n$-dimensional rational polytope $\Pc\!\subset\!\Rn$, and 
points $a_1,\ldots,a_k\!\in\!\Zn$, such that $\Pc$ is $n$-dimensional, 
$I$-cornered and compatible with $a\cup \cP$.} 
\item[Description:]{
\begin{enumerate}
\item{For all $i\!\in\![1..k]$ and $j\!\in\![1..n]$, define
$m_{ij}\!:=\!\min \limits \{e_j \; | \; (e_1,\ldots,e_n)\!\in\!P_i\}$, and 
let $m_1,\ldots,m_k$ be the rows of the matrix $[m_{ij}]$.} 
\item{For each $i\!\in\![1..k]$ let $a_i\!\in\!P_i\cap\Zn\cap (m_i+\lin(I))$ or
set $a_i\!:=\!m_i$ if $P_i\!\cap\!(m_i+\lin(I))\!=\!\emptyset$.} 
\item{Define $Q\!:=\!\sum^k_{i=1} \conv(\{a_i\}\!\cup\!P_i)$. If 
$\dim Q\!<\!n$ then set $Q\!:=\!Q+\conv(\bO\cup \cB)$, where $\cB$ is a 
generic set of $n-\dim Q$ rational points in the nonnegative orthant (so that 
$\dim Q\!=\!n$). }
\item{For all $i\!\in\![1..k]$, let $\eps_i\!:=\!\frac{1}{2}\min 
\{ {\mathrm v}_i \}$, where the minimum ranges over all vertices 
${\mathrm v}\!=\!({\mathrm v}_1,\ldots,
{\mathrm v}_n)$ of $Q$ incident to (but not lying in)
$(\hyp(i)+\sum^n_{j=1} a_j)\cap Q$. } 
\item{Define $\Pc\!:=\!((\eps_1,\ldots,\eps_n)+\sigma^\vee_I)\cap 
Q$.\footnote{It is useful to note that the dual cone 
$\sigma^\vee_I$ is precisely $\bigcap \limits_{j\in I^c} 
\{(y_1,\ldots,y_n)\!\in\!\Rn \; | \; y_j\geq 0\}$.} } 
\end{enumerate} }
\end{itemize}
\end{alg}
{}From the last step of our construction it is easily verified that $\Pc$
is $I$-cornered. Also, since $Q$ is already $n$-dimensional and compatible with
$a\cup\cP$, it is clear that our choice of $(\eps_1,\ldots,\eps_n)$ keeps 
$\Pc$ $n$-dimensional and compatible with $a\cup\cP$. Thus 
intersecting a translate of $\sigma^\vee_I$ with $Q$ in step (5) is somewhat 
reminiscent of refining the fan of a polytope by, quoting 
\cite[pg.\ 190]{gkz94}, ``cutting out (as with a knife)$\ldots$ any face of 
codimension at least 2.'' 

As one may have
already guessed, $\Pc$ is especially useful for root counting in $\Kni$ and
the points $a_1,\ldots,a_k$ will also be quite important. As a warm-up, the
following lemma is easily verified from theorem \ref{thm:orbit}, definition
\ref{dfn:impo2}, and proposition \ref{prop:corn}. 
\begin{lemma}
\label{lemma:embed} 
Following the notation of definition \ref{dfn:intersect}, assume further 
that $P_1,\ldots,P_k$ all lie in the nonnegative orthant of $\Rn$. Fix 
$I\!\subseteq\![1..n]$ and define $a_1,\ldots,a_k$ and $\Pc$ via algorithm 
\ref{alg:corner}. Then $\cP$ cornered
$\Longrightarrow (\Kni)\cap Z(F)\!=\!(\Kni)\cap
\cD_\Pc(F,a\cup\cP)$ as schemes. Furthermore, if $\newt(f_i)\!=\!P_i$ for 
all $i$ as well, then 
$\overline{(\Kni)\cap Z(F)}\!=\!\cD_\Pc(F,a\cup\cP)$. \qed 
\end{lemma}
\noindent
We emphasize that $\Kni$ is {\bf not} always naturally embedded in $\cT_P$, 
hence our need for $\Pc$. Thus, under certain assumptions, the above lemma 
allows us to embed an affine hypersurface into a toric divisor. In fact, we can 
do even better: We are now in a position to apply our framework to proving the 
Affine Point Theorem II.

\noindent
{\bf Proof of the Affine Point Theorem II:} Focusing on the first part of
the theorem, the case $\cM(a\cup E)\!=\!0$ is easiest to prove so we
dispose of it first: By the author's Affine Point Theorem I
\cite{rojaswang}, we obtain that a polynomial system with support contained
in $E$ can have {\em no}\/ isolated roots in $\Kni$. Since $E$ is
$(\Kni)$-nice by assumption, we are done.

So let us now assume that $\cM(a\cup E)\!>\!0$. Set 
$\cP\!:=\!(\conv(E_1),\ldots,\conv(E_n))$ and, applying algorithm 1, define  
$\cD\!:=\!\cD_\Pc(F,a\cup \cP)$. We will need the following important fact: 
\[ \star: \ O_w\cap\cD\!=\!\emptyset \mathrm{ \ for \ all \ } w\!\in\!\Rni
\Longleftrightarrow [\cD\!\subset\!\Kni \mathrm{ \ and \ } \dim
\cD\!\leq\!0]. \] 
\noindent
That the left-hand side is equivalent to $\cD\!\subset\!\Kni$ follows
easily from theorem \ref{thm:orbit} and proposition \ref{prop:corn}.
Furthermore, it follows easily from part (2) of the Antipodality Theorem
and proposition \ref{prop:corn} that $\dim \cD\!>\!0 \Longrightarrow
O_w\cap\cD\!\neq\!\emptyset$ for some $w\!\in\!\Rni$. So $\star$ is true.

Now note that the left-hand side of $\star$ is generically true by
proposition \ref{prop:over} and corollary \ref{cor:cool}. So by theorem
\ref{thm:toric}, lemma \ref{lemma:embed}, and $\star$, all but the last
sentence of the Affine Point Theorem II is now verified. Note 
also that we may drop the assumption that $E$ be $(\Kni)$-nice, 
as long as we also count {\em embedded}\/ zero-dimensional components.  

To prove the final part, first note the following identity of weighted set
unions: 
\[ \dagger: \ O_J = (\Kn\!\setminus\!\hyp(J)) \setminus
\left( \bigcup\limits_{|J'\setminus J|=1} (\Knjp) \right) \cup \left(
\bigcup\limits_{|J'\setminus J|=2} (\Knjp) \right) \setminus \cdots \]
which terminates in the appropriate union or set difference according as
$n-|J|$ is even or odd. This follows easily from the principle of
inclusion-exclusion \cite{gkp} since, for any $\vartheta\!\subseteq\![1..n]$, 
$\Kn\!\setminus\!\hyp(\vartheta)$ is precisely the disjoint
union $\coprod_{\vartheta'\supseteq \vartheta} O_{\vartheta'}$. The key to 
proving our alternating mixed volume formula is then to simply find an 
intersection theoretic analogue of $\dagger$.

To do this we must work in a new {\em lifted}\/ compactification depending on 
$J$.  So let $P(J')$ denote the $\Pc$ corresponding to the $I\!=\!J'$ 
case of algorithm 1 and define $\wcT$ to be the toric compactification
corresponding to $\wP\!:=\!\sum_{J'\supseteq J} P(J')$. Also let
$a_1(J'),\ldots,a_n(J')$ respectively denote the integral points
$a_1,\ldots,a_n$ from the $I\!=\!J'$ case of algorithm 1. By example 
\ref{ex:compat}, $\wP$ is compatible with $a(J')\cup\cP$ for all 
$J'\!\supseteq\!J$, so 
define $\wcD(J')\!:=\!\cD_\wP(F,a(J')\cup\cP)$. Then, by \cite[Chap.\
2.4]{tfulton} and our construction, there is a proper morphism 
$\varphi : \wcT \twoheadrightarrow \cT_{P(J)}$ with no fibers of 
infinite cardinality. Similar to proposition \ref{prop:corn} and lemma 
\ref{lemma:embed}, it is also easily checked that 
$O_J\cap\varphi(\wcD(J))\!=\!O_J\cap Z(F)$.  

More importantly, it is easily verified from definition \ref{dfn:impo2} 
and expanding in $\chow(\wcT)$ that (as cycles) $\wcD(J'')$ is a 
{\em summand}\/ of $\wcD(J')$ for all
$J''\!\supseteq\!J'\!\supseteq\!J$. Also, it immediately follows from
theorem \ref{thm:toric} that $\deg \wcD(J')\!=\!\cM_{J'}$ for all
$J'\!\supseteq\!J$. So by inclusion-exclusion once again, we have the
following equality of cycles: \[ \varphi^{-1}(O_J)\cap \wcD(J) = \wcD(J) - 
\left( \sum\limits_{|J'\setminus J|=1}
\wcD(J') \right)
+ \left( \sum\limits_{|J'\setminus J|=2}
\wcD(J') \right)
- \cdots \]
provided $\varphi^{-1}(O_J)\cap \wcD(J)$ is zero-dimensional or 
empty.\footnote{We should remark that the left-hand intersection 
is {\em set}-theoretic, and {\em not}\/ a Chow product.} 
So then $\deg(\varphi^{-1}(O_J)\cap\wcD(J))$ is precisely our 
alternating mixed volume formula. Note that $\wcD(J)$ generically 
has exactly $\deg(\varphi^{-1}(O_J)\cap\wcD(J))$ points 
(counting multiplicities) in $\varphi^{-1}(O_J)$, by the portion of 
the Affine Point Theorem II that we've already proved and since 
$E$ is $O_J$-nice. The last cycle class degree is also precisely 
$\deg(O_J\cap\varphi(\wcD(J)))$ \cite[Example 7.1.9]{ifulton2}, so we are 
done. \qed

Of course, the assumption that $E$ be cornered is quite restrictive. We
relax this assumption in the following section by refining lemma
\ref{lemma:embed}, and then Main Theorem 1 follows easily by explicitly
expanding a different intersection product in the Chow ring of $\cT_\Pc$.

\subsection{Chow Rings and Main Theorem 1} 
\label{sub:chow}
Our second and final trick for applying special $\cT_P$'s to affine root
counting is a bit more abstract. Whereas our first trick (``cornering'') 
consisted of a convex geometric construction, the construction we give now 
amends a difficulty with the divisors $\cD_P(f,Q)$ we used earlier. In 
particular, for noncornered $(P_1,\ldots,P_n)$, it is possible that 
$\overline{(\Kni)\cap Z(f_i)}$ and $\cD_\Pc(f_i,\conv(\{a_i\}\cup P_i))$ 
differ in the coefficients corresponding to the coordinate hyperplanes. This 
is remedied by the following definition and lemma.

\begin{dfn}
\label{dfn:shift}
Following the notation of Main Theorem 1, definition \ref{dfn:toric}, and 
lemma \ref{lemma:embed}, define
$\cX_j\!:=\!V_{\hat{e}_j} \subset \cT_\Pc$ for any
$j\!\in\!I^c$. Also, set $\Ds(0,P_i)\!:=\!\cT_\Pc$ and, if 
$\supp(f_i)\!\neq\!\emptyset$, define $\Ds(f_i,P_i)\!:=\!\cD_\Pc(f_i,
\conv(\{a_i\}\cup P_i))+\sum_{j\in I^c} 
m_{ij}\cX_j \in \chow(\cT_\Pc)$. Finally, let $\Ds(F,\cP)\!:=\!\bigcap^k_{i=1} 
\Ds(f_i,P_i)\in\chow(\cT_\Pc)$. The {\em roots of $f_i$ within $\cT_\Pc$} 
are then, formally, $\Ds(f_i,\newt(f_i))$. 
\end{dfn}
\begin{lemma}
\label{lemma:shift}
Following the notation of definition \ref{dfn:shift}, $(\Kni)\cap
Z(F)=(\Kni)\cap\Ds(F,\cP)$ as schemes. Furthermore, if 
$\newt(f_i)\!=\!P_i$ then $\overline{(\Kni)\cap Z(f_i)}=\Ds(f_i,P_i)$, within 
$\cT_\Pc$. \qed
\end{lemma}
\noindent
In particular, for any $j\!\in\!I^c$, $\cX_j$ is the
closure of the hyperplane $\{ x \; | \; x_j=0\}\cap(\Kni)$ in $\cT_{\Pc}$.
Keeping this in mind, the proof of the lemma is then straightforward from
theorem \ref{thm:orbit}, proposition \ref{prop:corn}, and lemma
\ref{lemma:embed}. So by ``shifting'' our toric divisors, we now at last have 
a completely general way of embedding an affine hypersurface into a toric 
compactification. As an application, we will prove Main Theorem 1.

\noindent
{\bf Proof of Main Theorem 1:} Set $\cP\!:=\!(\conv(E_1),\ldots,\conv(E_n))$. 
We will first prove the case $W\!=\!\Kni$ and, to do so, it will clearly 
suffice to demonstrate the following two statements:
\begin{itemize}
\item[A$_{\deg}$: ]{$\cN_K(E;\Kni)\!=\!\deg\Ds(F,\cP)$. }
\item[A$_{\sm}$: ]{$\deg \Ds(F,\cP)$ is precisely the double summation stated 
in Main Theorem 1.}
\end{itemize}
Consider also the following auxiliary statement: 
\begin{itemize}
\item[A$_{\gen}$: ]{for fixed $E$ and generic $\cC_E$, $\Ds(F,\cP)$ is 
zero-dimensional and supported entirely within $\Kni$.} 
\end{itemize}
To prove A$_{\deg}$ and A$_{\sm}$, we will actually 
first prove $($A$_{\deg})\wedge($A$_{\gen})$ by induction on $n$, 
and then A$_{\sm}$ will follow easily. 

First note that by the definition of $\Ds(\cdot)$ we may formally expand
$\Ds(F,\cP)$ in $\chow(\cT_\Pc)$ as a polynomial in the
$\cX_j$. More explicitly, 
\[ \Ds(F,\cP)=\bigcap\limits^n_{i=1}\left(\cD_\Pc(f_i,\conv(\{a_i\}\cup P_i)) + 
\sum\limits_{j\in I^c} m_{ij}\cX_j\right)= \]
\[\sum \limits_{[1..n]\supseteq J\supseteq I} \; 
\sum \limits_{\rho:J^c \hookrightarrow [1..n]} \left[
\left( \prod \limits_{j\in J^c} m_{j\rho(j)} \right)
\left(\bigcap \limits_{i\in J}\cD_\Pc(f_i,\conv(\{a_i\}\cup P_i))\right)
\cap\left( \bigcap \limits_{j\in \rho(J^c)}\cX_j \right)\right] \]
This is where the shape of our asserted formula comes from. Note 
that $j\!\in\!I \Longrightarrow \cX_j\cap(\Kni)\!=\!\emptyset$, thus allowing 
the slight simplification of the outer summation.  

Now note that $\bigcap_{j\in \rho(J^c)} \cX_j$ is itself isomorphic 
to the toric variety corresponding to a face $P_{(J,\rho)}$ of $\Pc$,  
\`a la theorem \ref{thm:orbit}. In particular, letting 
$\cP_{(J,\rho)}\!:=\!((P_j-m_j)\cap\lin(\rho(J^c)^c) \; | \; 
j\!\in\!J)$, it easy to see that $P_{(J,\rho)}$ can occur as the output 
of the $(k,n,\cP,I)\!\rightsquigarrow\!(|J|,|J|,\cP_{(J,\rho)},
\rho(J^c)^c\cap I)$ case of algorithm 1. 

Let $F_{(J,\rho)}$ be the polynomial system obtained by 
setting the variables $\{x_j \; | \; j\!\in\!\rho(J^c)\}$ to
$0$ in the $|J|$-tuple $(x^{-m_{i1}}_1\cdots x^{-m_{in}}_nf_i \; | \; 
i\!\in\!J)$. Also note that 
$\bigcap_{j\in \vartheta}\hyp(j)\!=\!\lin(\vartheta^c)$ for any 
$\vartheta\!\subseteq\![1..n]$. We may then say that 
\[ \left( \bigcap 
\limits_{i\in J}\cD_\Pc(f_i,\conv(\{a_i\}\cup P_i))\right) 
\cap\left(\bigcap \limits_{j\in \rho(J^c)}\cX_j\right) \cong 
\Ds(F_{(J,\rho)},\cP_{(J,\rho)}) \] 
where the underlying compactification for the right-hand cycle is 
$\cT_{P_{(J,\rho)}}$. This last identity follows from definitions 
\ref{dfn:impo2} and \ref{dfn:shift}, and our preceding observations. 

Note that $E_{([1..n],\cdot)}$ is cornered. So then the proof of the Affine 
Point Theorem II (and definition \ref{dfn:impo2}) immediately implies that 
$\deg\cD_\Pc(F,a\cup\cP)\!=\!\cM(a\cup\cP)\!=\!\cN_K(E_{([1..n],\cdot)};\Kni)$ 
and, generically, $\cD_\Pc(F,a\cup\cP)$ is zero-dimensional and 
supported entirely within $\Kni$.  As for the remaining intersection terms 
with $J\!\neq\![1..n]$, our induction hypothesis (with $n\!=\!|J|$) implies 
that
\[\cN_K\left(E_J;\lin(\rho(J^c)^c)\cap(\Kni)\right)\!=\!\deg\Ds(F_{(J,\rho)},
\cP_{(J,\rho)}).\] 
Furthermore, our induction hypothesis also implies that, generically, 
$\Ds(F_{(J,\rho)},\cP_{(J,\rho)})$ is zero-dimensional and supported entirely 
within $\lin(\rho(J^c)^c)\cap(\Kni)$. 

Now note that our Chow expansion also immediately implies that 
\[
\supp(\Ds(F,\cP))\!=\!\bigcup_{I\subseteq J\subseteq[1..n]} \;  
\bigcup_{\rho:J^c\hookrightarrow [1..n]} \supp(\Ds(F_{(J,\rho)},\cP_{(J,\rho)})), \] 
modulo some isomorphisms fixing $\Kni$. Since a finite conjunction of generic
conditions is again a generic condition, we thus 
arrive at A$_{\gen}$.  

Recall that lemma \ref{lemma:shift} states 
that $(\Kni)\cap Z(F)$ is naturally embedded in $\Ds(F,\cP)$. Thus 
$\cN_K(E;\Kni)\!\leq\!\deg\Ds(F,\cP)$ and, by A$_{\gen}$, we arrive at 
A$_{\deg}$. Noting that the $n\!=\!1$ case of $($A$_{\deg})\wedge($A$_{\gen})$ 
is true simply via the fundamental theorem of algebra over $K$, our 
induction is complete. 

Finally, A$_{\sm}$ follows simply by taking degrees of 
both sides of our Chow expansion. Note also that our embedding, along with 
A$_{\deg}$, implies that $\cN_K(E;\Kni)$ is indeed the maximal number of 
isolated roots. So the case $W\!=\!\Kni$ is proved. 

The general case then follows easily from inclusion-exclusion, much like
our proof of the Affine Point Theorem II. This method goes through because
our asserted formula is additive with respect to disjoint unions in $W$, and 
already true for $W\!=\!\Kni$. \qed
\begin{rem}
Since our proof of Main Theorem 1 computes $\cN_K(E;W)$ as the degree 
of an algebraic cycle, remark \ref{rem:disting} is just a 
straightforward abstract extension of our preceding proof.
\end{rem}
\begin{rem}
\label{rem:simple}
The double summation of Main Theorem 1 can of course simplify considerably 
when $W$ is smaller than $\Kn$. For instance, there is only 1 term when 
$W\!=\!\Ksn$.  Also, it is a simple combinatorial exercise to show that the 
double summation of Main Theorem 1 has at most $\prod^n_{i=1} 
(|\supp(m_i)|+1)$ terms. Equality occurs, for example, when $W\!=\!\Kn$. 
Note also that by our recursive formula, the matrix $[m_{ij}]$ can be assumed to
have at most $1$ nonzero entry per column if 
$E$ is $W$-nice. So, by the fact that $(a+1)(b+1)\!\geq\!a+b+1$ for positive 
integers, we also obtain that our double summation has at most $2^n$ terms. 
Furthermore, this maximum is attained iff $[m_{ij}]$ has the same support as a 
permutation matrix.
\end{rem}
\begin{rem}
If one would like a formula closer to the number of {\em distinct} roots in
$\Kni$, a useful trick is the following: Use Main Theorem 1,
but replacing $[m_{ij}]$ with the 0-1 matrix having the same support. From
our last proof, it is easy to see that this new formula has the effect of
(generically) counting isolated roots lying on $\bigcup \limits_{J\supseteq 
I, \ |J|=n-1} O_J$ {\em without} multiplicity. Thus, if one continues 
to propogate this trick throughout the recursion of Main Theorem 1, 
we can count (omitting multiplicities due to inseparability degree) the 
generic number of distinct roots of $F$ in $W$. 
This is important because for many $E$, a sparse system with support contained 
in $E$ {\em always}\/ has roots of multiplicity $>\!1$ lying in 
$\Kn\!\setminus\!\Ksn$. 
\end{rem}

A useful corollary of our proof of Main Theorem 1 is the following concise
generalization of theorem \ref{thm:toric}. 
\begin{cor}
\label{cor:new}
Following the notation of lemma \ref{lemma:shift}, 
assume further that $k\!=\!n$ and $\cP\!=\!(\conv(E_1),\ldots,\\ 
\conv(E_n))$. Then $\cN_K(E;\Kni)\!=\!\deg \Ds(F,\cP)$. Furthermore, if both 
$\cN_K(E;\Kni)$ and $\dim \Ds(F,\cP)$ are positive, then $\Ds(F,\cP)$ has 
strictly less than $\cN_K(E;\Kni)$ zero-dimensional components, counting 
multiplicities. 
\end{cor}
\noindent
{\bf Proof:} The first portion follows immediately from our proof of the 
$W\!=\!\Kni$ case of Main Theorem 1. As for the remaining portion,  
by the Chow expansion from our last proof,  
it suffices to prove the cornered case and then simply mimic the earlier 
descent by induction. Since the cornered case of our present corollary is 
already contained in the $(\cP,P)\!\rightsquigarrow\!(a\cup\cP,P^\llcorner)$ 
case of theorem \ref{thm:toric}, we are done. \qed

So theorem \ref{thm:toric} is just the $I\!=\![1..n]$ case of 
collary \ref{cor:new}. In intersection theoretic terms, the above
result establishes the {\em numerical positivity}\/ \cite{ifulton2} of 
any positive-dimensional component of the new shifted cycle $\Ds(F,\cP)$. 
This will allow us to derive precise algebraic conditions for what ``generic'' 
means in the context of affine root counting.

Corollary \ref{cor:stable} then follows easily from Main Theorem 1 as follows: 

\noindent
{\bf Proof of Corollary \ref{cor:stable}:} Although Huber and Sturmfels did
not explicitly mention intersection multiplicities in \cite{hsaff}, an
examination of their proof of the stable mixed volume formula shows that 
multiplicities were at least included implicitly. In particular, we may
safely assume that the first portion of Corollary \ref{cor:stable} is true
for $K\!=\!\C$. The remaining portion (for $K\!=\!\C$) is already implicit
in Huber and Sturmfels' proof\footnote{Note that our $O_J$ is actually 
$O_{\{1,\ldots,n\}\setminus J}$ in the notation of \cite{hsaff}.} 
of the stable mixed volume formula, so we may safely assume that {\em all}\/ 
of Corollary \ref{cor:stable} is true for $K\!=\!\C$.

Generalizing to arbitrary algebraically closed $K$ is then almost trivial: 
The right-hand sides (of both asserted formulae) are clearly independent of
$K$. By Main Theorem 1 and the Affine Point Theorem II, the left-hand sides
are also independent of $K$, provided $K$ is algebraically closed. Since
both formulae are already true for $K\!=\!\C$, we are done. \qed

Similar to remark \ref{poly}, a more elementary (but longer) proof of
Corollary \ref{cor:stable} can be derived by generalizing Huber and Sturmfels' 
proof of their stable mixed formula. 

\subsection{Sparse Resultants, Roots at Infinity, and Main Theorems 2 and 
3} 
\label{sub:algcond} 
We conclude with an analysis of conditions under which 
our (global) generic root counts are exact. The conditions 
we give can be split into two types: algebraic and combinatorial. 
In the combinatorial case our conditions are always both sufficient and 
necessary, while in the algebraic case our criteria are always sufficient 
but fail to be necessary for certain systems which generically have no roots. 
However, we fully classify the cases where our algebraic criteria are 
necessary. These results will rely on the following 
technical result relating our shifted toric divisors with toric infinity. 
\begin{lemma}
\label{lemma:split}
Following the notation of the proof of Main 
Theorem 1, let $\cD_{[1..n]}$ ($\in\!\chow(\cT_\Pc)$) be the $J\!=\![1..n]$ 
term of the Chow expansion of $\Ds(F,\cP)$. Then the following conditions 
imply that $F$ has exactly $\cN_K(E;\Kni)$ roots, counting multiplicities, in 
$\Kni$: 
\begin{itemize}
\item[(a)]{$O_w\cap\cD_{[1..n]}\!=\!\emptyset$ 
for all $w\!\in\!\Rni$, and }    
\item[(b)]{if $n\!>\!1$ then for all $J\!\subsetneqq\![1..n]$ containing
$I$, and all injections $\rho : J^c \hookrightarrow [1..n]$ such that
$\rho(J^c)\cap I\!=\!\emptyset$ and $\prod_{j\in J^c} m_{j\rho(j)}>0$, 
\[\cN_K\left(E_{(J,\rho)};\lin(\rho(J^c)^c)\cap(\Kni);\cC_E\right)=
\cN_K\left(E_{(J,\rho)};\lin(\rho(J^c)^c)\cap(\Kni)\right).\]}
\end{itemize}
Furthermore, the converse implication holds as well if 
$\cN_K(E_{([1..n],\cdot)};\Kni)\!>\!0$. In particular, (a) and (b) 
together imply that the zero set of $F$ in $\Kni$ is zero-dimensional or 
empty. 
\end{lemma} 
\noindent
{\bf Proof of the Lemma:} 
Note that $\supp(\Ds(F,\cP))\!=\!\bigcup \supp(\Ds(F_{(J,\rho)},
\cP_{(J,\rho)}))$ where the union ranges over $([1..n],\cdot)$ and all pairs 
$(J,\rho)$ described above. This follows immediately from our Chow 
expansion from the proof of Main Theorem 1, and the fact that the 
terms with $\rho(J^c)\cap I$ nonempty or $\prod_{j\in J^c} m_{j\rho(j)}$ 
zero simply aren't there (by definition \ref{dfn:shift}). Recall also that 
$E_{([1..n],\cdot)}$ is cornered. So it suffices to prove the cornered case 
and descend by induction, just as we did in the proofs of Main Theorem 1 and 
corollary \ref{cor:new}. But the cornered case, minus 
the partial converse, is already contained in assertion $\star$ from our 
proof of the Affine Point Theorem II. One also observes that the 
$(\Longrightarrow)$ portion of $\star$ continues to hold even when 
$\cN_K(E;\Kni)\!=\!0$. So we are done. \qed 
\begin{rem}
\label{rem:fancy}
Note that the converse of the main assertion of lemma \ref{lemma:split} 
can fail if $\cN_K(E_{([1..n],\cdot)};\Kn\!\setminus\\
\hyp(I))\!=\!0$: Consider the 
polynomial system $F\!=\!(1+x,1+x,(1+x)(y+z)+1)$ where $E\!:=\!\supp(F)$ 
and note that condition (b) is violated when $w\!=\!(0,-1,-1)$. 
\end{rem}
\begin{rem}
\label{rem:fancier}
However, the converse doesn't always fail if 
$\cN_K(E_{([1..n],\cdot)};\Kni)\!=\!0$: For instance, the converse 
{\em always} holds for $n\!=\!2$ when $E\!=\!\supp(F)$. More generally,   
for any $n\!>\!2$, setting $f_1\!:=\!1$ and $E_1\!:=\!\{\bO\}$ gives 
an entire family of examples. Basically, when $\cN_K(E;\Kni)\!=\!0$, 
the converse of the main assertion of lemma \ref{lemma:split} fails precisely 
when $E$ is sufficiently complicated to allow specializations of $\cC_E$ where 
$F$ has roots at toric infinity while having {\em none} within $\Kni$. 
\end{rem}

We are now ready to prove Main Theorem 3.

\noindent 
{\bf Proof of Main Theorem 3:} We will first dispose of case (1) which 
is the easiest. Recall that $\Kni\!=\!\coprod_{J\supseteq I}
O_J$ and that a finite conjunction of generic conditions is again a
generic condition. Since $E$ is null for $\Kni$ (and thus for every 
$O_J$ with $J\!\supseteq\!I$) it suffices to show that 
our condition from case (1) is equivalent to $D\cap\lin(J)$ $O_J$-counting 
$E\cap\lin(J)$ for all $J\!\supseteq\!I$. This, in essence, is the statement 
of Lemma 3 of \cite{convexapp}. So case (1) is complete. 

As for case (2), note that $E^w$ depends only on the face $S^w$. So by 
corollary \ref{cor:cool}, the same is true of $O_w\cap\cD_{[1..n]}$. 
Then by lemma \ref{lemma:split}, the definition of 
$W$-counting, and since any finite conjunction of generic conditions is again 
a generic condition, we need only prove the case where $E$ is cornered 
and then descend by induction just as in three of our last four proofs. For 
$I\!=\!\emptyset$, the cornered case is just case (2) of Theorem 7 of 
\cite{rojaswang}. Applying algorithm 1, generalizing the proof 
there to arbitrary $I$ is simple. (In fact, the proof of 
\cite[Theorem 7]{rojaswang} already contains what is essentially the 
$I\!=\!\emptyset$ case of algorithm 1.) So we are done. \qed 

We now recall the {\em sparse resultant}\/ (also known as the 
{\em $(\cA_1,\ldots,\cA_k)$-resultant}, {\em mixed}\/ resultant, {\em 
Newton}\/ resultant, or {\em toric}\/ resultant), which is 
an extremely important operator on overdetermined polynomial systems. It is 
defined for any $k\!\times\!n$ {\bf indeterminate} polynomial system
$F$ with support $E$, provided that all the $E_i$ can be translated into a
common $(k\!-\!1)$-dimensional subspace of $\Rn$. Since we can always 
identify such a subspace with a rational hyperplane in $\R^k$, we will   
consider only the case of $n\!\times\!(n-1)$ systems and 
monomial transformations (involving an extra variable) of such systems. 

More explicitly, suppose $E$ is an $n$-tuple of nonempty finite subsets 
of $\Zn$ which can be translated into a common $(n\!-\!1)$-plane in $\Rn$. 
Then the sparse resultant, with respect to $E$, will be a (homogeneous) 
polynomial $\res_E(\cdot)$ in the coefficients $\cC_E$ satisfying the 
following property: If $\cC\!\in\!K^{|E|}$ and $F|_{\cC_E=\cC}$ has a root in 
$\Ksn$, then $\res_E(\cC)\!=\!0$. For 
fixed $E$, the polynomial $\res_E(\cdot)$ can then be {\bf defined} (up to a 
nonzero scalar multiple) as the unique polynomial in $\cC_E$ of least total 
degree satisfying this last property. The computation of 
$\res_E(\cdot)$ is a deep subject and we refer the reader to
\cite{gkz90,chowprod,jandi1,sz94, combiresult,gkz94,isawres,resvan} for
further background on sparse resultants.

For convenience, we will use $\res_E(F)$ in place of $\res_E(\cC)$ whenever
the coefficients of $F$ have been specialized to some $\cC\!\in\!K^{|E|}$.
We also point out the following important fact: $\res_E(F)\!=\!0$ does 
{\bf not} necessarily imply that $F$ has a root in $\Ksn$. The correct
statement, at least for initial term systems, is the following.
\begin{thm}
\cite{resvan}
\label{thm:resorbit}
Following the notation of theorem \ref{thm:toric}, suppose 
$\cP\!=\!(\conv(E_i) \; | \; i\!\in\![1..n])$ and $w$ is an inner
facet normal of $P$. Then $\res_{E^w}(F)\!=\!0
\Longleftrightarrow V_w\!\cap\!\cD_P(F,\cP)\!\neq\!\emptyset$. \qed 
\end{thm}
\begin{rem}
Since $\init_{w,E}(F)$ involves a subset of the coefficients $\cC_E$, we've
further simplified notation by writing $\res_{E^w}(F)$ in place of
$\res_{E^w}(\init_{w,E}(F))$. 
\end{rem}

Our algebraic condition for $F$ to have generically many roots 
is based on the above useful property of the sparse resultant.

\noindent
{\bf Proof of Main Theorem 2:} Recall from the proof of Main Theorem 
3 that $E^w$ and $O_w\cap\cD_{[1..n]}$ depend only on the face $S^w$.  
So then, by theorem \ref{thm:orbit}, $V_w\cap\cD_{[1..n]}$ also depends 
only on $S^w$. To conclude, by theorems \ref{thm:orbit} and \ref{thm:resorbit}, 
it is clear that conditions (a) and (b) of lemma \ref{lemma:split} are 
respectively equivalent to conditions (a$_2$) and (b$_2$) of Main Theorem 2. 
So we are done. \qed 
\begin{rem}
In practice, we would not actually construct the Minkowski 
sum $S$ stated in our last two main theorems. Instead, we would actually  
work with the individual normal fans of $\conv(E_1),\ldots,\\
\conv(E_n)$ and 
dynamically update the smallest common refinement necessary for our search 
space. The deeper stages in the recursion would then simply consist of 
slicing cones in the search fans by appropriate coordinate subspaces.  
\end{rem} 

The case $(K,I)\!=\!(\C,[1..n])$ of Main Theorem 2 was independently 
discovered and presented in \cite[Theorem 6.1]{polyhomo}. However, the 
statement there is false in the case where the mixed volume is zero: Simply 
consider the counter-example from remark \ref{rem:fancy}. Moreover, parallel 
to lemma \ref{lemma:split}, the converse of the main assertion of Main Theorem 
2 doesn't depend completely on the positivity of $\cN_K(E;\Kni)$: The examples 
given in remarks \ref{rem:fancy} and \ref{rem:fancier} also work here in an 
analogous way. Sharper computational conditions for exactness in the cases  
$\cN_K(E;\Kni)\!=\!0$ will be addressed in future work. 
\begin{rem} 
Following the notation of definition \ref{big}, we see that  
Main Theorem 2 allows us to express $\Delta$ as a union of 
hypersurfaces, via condition (a$_2$) and the recursion of condition 
(b$_2$). However, the {\em minimal} delta (containing 
all $\cC$ such that $F|_{\cC_E=\cC}$ does {\em not} have the 
generic number of roots) need not be a hypersurface or even an 
algebraic variety. In general, this ``generic counting discriminant'' is a 
{\em constructible} variety of codimension $>\!1$. However, we 
can at least explicitly compute the codimension by combining 
Main Theorem 2 with Theorem 1.3 of \cite{combiresult}. 
\end{rem}

``Sparse'' techniques have recently been applied quite succesfully to {\em 
solving}\/ many polynomial systems occuring in industrial problems
\cite{mobiorobo,emiphd,isawres,dynamiclift}. Software implementations of
resultant-based algorithms are also discussed in almost all of these
papers. Thus Main Theorem 2 presents another potentially useful 
application of the sparse resultant. 
 
\section{Acknowledgements}
\label{sec:acknowledgements}
The author would like to thank Marie-Francoise Roy for her kind 
support and hospitality during his visit to IRMAR in Rennes. 
He is also grateful to Paco Santos and Jie Tai Yu for shooting down 
his erroneous conjectures, and to Birk Huber, T. Y. Li, Bernd
Sturmfels, and Robert M. Williams for many valuable discussions. The author 
also thanks the referees from MEGA '96 for their valuable suggestions.  

\bibliographystyle{amsalpha}

\section*{Appendix: Niceness and Genericity}
Here we briefly recount some earlier results on $W$-counting and some
related concepts. Some of the material below is covered at greater length
in \cite{convexapp} (for the case $I\!=\![1..n]$) and \cite{rojaswang} (for
the case $I\!=\!\emptyset$ and $E$ cornered). The paper \cite{combiresult}
is also a useful reference but deals more with the sparse resultant than
with root counting. The results below form the basis for our combinatorial
conditions for when a ``partially'' generic polynomial system has
generically many isolated roots in a given union of orbits. 

Recall that the {\em dimension}\/ of any $B\!\subseteq\!\Rn$, $\dim B$, is
the dimension of the smallest subspace of $\Rn$ containing a translate of
$B$. The following two definitions are fundamental to our development.
\begin{dfn}
Suppose $C\!:=\!(C_1,\ldots,C_n)$ is an
$n$-tuple of polytopes in $\Rn$ {\em or} an $n$-tuple of finite subsets of
$\Rn$. We will allow any $C_i$ to be empty and say that a nonempty subset
$J\!\subseteq\![1..n]$ is {\em essential} for $C$ (or $C$ {\em has
essential subset $J$}) $\Longleftrightarrow$ (0) $\supp(C)\!\supseteq\!J$,
(1) $\dim(\sum_{j\in J} C_j)=|J|-1$, and (2) $\dim(\sum_{j\in J'} C_j)\geq
|J'|$ for all nonempty {\em proper} $J'\!\subsetneqq\! J$.
\end{dfn}
\begin{dfn}
We say that $C$ {\em has an almost essential subset $J$}
$\Longleftrightarrow$ (0) $\supp(C)\!\supseteq\!J$, (1) $\dim(\sum_{j\in J}\\
C_j)=|J|$, and (2) $\dim(\sum_{j\in J'} C_j)\geq |J'|$ for all nonempty 
$J'\!\subseteq\! J$. Also, $\emptyset$ is defined to be almost essential for 
$C$ iff $\supp(C)\!=\!\emptyset$. 
\end{dfn}

Equivalently, $J$ is essential for $C \Longleftrightarrow$ the
$|J|$-dimensional mixed volume of $(C_j \; | \; j\!\in\!J)$ is $0$ and no
smaller subset of $J$ has this property. The following figure shows some
simple examples of essential subsets for $C$, for various $C$ in the case
$n\!=\!2$. It also worth noting that $J\cup\{j\}$ is essential for $C 
\Longrightarrow J$ is almost essential for $C$, provided
$|J\cup\{j\}|\!>\!|J|\!>\!0$. 

Making use of the fact that $\Kni$ is the disjoint union
$\coprod_{J\supseteq I} O_J$, our combinatorial results for
$(\Kni)$-counting and $(\Kni)$-niceness will follow easily upon
partitioning $\Kni$ into orbits. In particular, it will be useful to refine
$W$-niceness slightly as follows. 
\begin{dfn}
Suppose $E$ is an $n$-tuple of finite subsets of $(\N\cup\{0\})^n$. 
We then call $E$ {\em null} for $W \Longleftrightarrow$ a generic polynomial 
system with support contained in $E$ has {\em no} roots in $W$. 
\end{dfn}

For any $J\!\subseteq\![1..n]$, define
$E\cap\lin(J)\!:=\!(E_1\cap\lin(J),\ldots,E_n\cap\lin(J))$. We may now
quote the following useful result.
\begin{lemma}
\label{lemma:start}
\cite[Corollary 2]{rojaswang}
Suppose $E$ is an $n$-tuple of finite subsets of $(\N\cup\{0\})^n$. Then
$E$ is nice for $O_J \Longleftrightarrow E\cap\lin(J)$ has an almost
essential subset of cardinality $|J|$ or an essential subset. In
particular, $E$ is null for $O_J \Longleftrightarrow E\cap\lin(J)$ has an
essential subset. \qed 
\end{lemma}

The following characterization of $(\Kni)$-niceness then follows almost
immediately.
\begin{lemma}
\label{lemma:wow}
An $n$-tuple $E$ of finite subsets of $(\N\cup\{0\})^n$ is nice for $\Kni
\Longleftrightarrow$ for all $J\!\supseteq\!I$, $E\cap\lin(J)$ has an
almost essential subset of cardinality $|J|$ or an essential subset. In
particular, $E$ is null for $\Kni \Longleftrightarrow$ for all
$J\!\supseteq\!I$, $E\!\cap\!\lin(J)$ has an essential subset. \qed
\end{lemma}
\noindent
The characterization of $W$-niceness for $W$ an arbitrary union of orbits 
is then completely analogous. 

\end{document}